\documentclass{article}
\usepackage{amsmath, epsfig,amssymb, multicol}
\newtheorem{lemma}{Lemma}
\newtheorem{theorem}{Theorem}

\newtheorem{corollary}{Corollary}

\newcommand{\C}{{\cal C}}

\newcommand{\F}{{\cal F}}
\renewcommand{\c}{{\bf c}}

\newcommand{\I}{{\cal I}}

\title{The Fractional Chromatic Number of Triangle-free Graphs
with $\Delta\leq 3$}
\author{Linyuan Lu
\thanks{University of South Carolina, Columbia, SC 29208,
({\tt lu@math.sc.edu}).
This author was supported in part by NSF grant
DMS 0701111 and DMS 1000475.}
  \and Xing Peng
\thanks{University of South Carolina, Columbia, SC 29208,
({\tt pengx@mailbox.sc.edu}).This author was supported in part by NSF grant
DMS 0701111 and DMS 1000475.}
}
\begin{document}
\maketitle
\begin{abstract}
Let $G$ be a  triangle-free graph with maximum degree at most 3.
Staton proved that the independence number of $G$ is at least
$\frac{5}{14}|V(G)|$. Heckman and Thomas conjectured that Staton's
result can be strengthened into a bound on the fractional chromatic
number of $G$, namely  $\chi_f(G)\leq \frac{14}{5}$. Recently,
Hatami and Zhu proved that $\chi_f(G) \leq 3 -\frac{3}{64}$. In this
paper, we prove $\chi_f(G) \leq 3- \frac{3}{43}$.
\end{abstract}
\section{Introduction}

This paper investigates the fractional chromatic number of a
triangle-free graph with maximum degree at most $3$. For a simple
(finite) graph $G$, the fractional chromatic number of  $G$ is the
linear programming relaxation of the chromatic number of $G$. Let
$\I(G)$ be the family of independent sets of $G$. A mapping $f\colon
\I(G)\to [0,1]$ is called an $r$-fractional coloring of $G$ if
$\sum_{S\in \I(G)}f(S)\leq r$ and $\sum_{v\in S, S\in \I(G)}f(S)\geq
1$ for each $v \in V(G)$. The {\it fractional chromatic number} $\chi_f(G)$ of $G$ is the
least $r$ for which $G$ has an $r$-fractional coloring.

Alternatively, the fractional chromatic number can also be defined
through multiple colorings.  A $b$-{\it fold coloring} of $G$
assigns a set of $b$ colors to each vertex such that any two
adjacent vertices receive disjoint sets of colors. We say a graph
$G$ is $a$:$b$-{\it
  colorable} if there is a $b$-fold coloring of $G$ in which
each color is drawn from a palette of $a$ colors. We refer to such a
coloring as an $a$:$b$-coloring. The $b$-{\it fold coloring number},
denoted as $\chi_b(G)$, is the smallest integer $a$ such that $G$
has a $a$:$b$-coloring. Note that $\chi_1(G)=\chi(G)$. It is known
\cite{su} that $\chi_b(G)$ (as a function of $b$) is sub-additive
and so the $\displaystyle\lim_{b\to\infty}\tfrac{\chi_b(G)}{b}$ always exists, which
turns out to be an alternative definition of $\chi_f(G)$. (Moreover,
$\chi_f(G)$ is a rational number and the limit can be replaced by
minimum.)

Let $\chi(G)$ be the chromatic number of $G$ and $\omega(G)$ be the
clique number of $G$. We have the following simple relation,
\begin{equation}
  \label{eq:1}
\omega(G)\leq \chi_f(G)\leq \chi(G).
\end{equation}

Now we consider a graph $G$ with maximum degree $\Delta(G)$ at most
three. If $G$ is not $K_4$, then  $G$ is $3$-colorable by Brooks'
theorem. If $G$ contains a triangle, then $\chi_f(G)\geq
\omega(G)=3$. Equation (\ref{eq:1}) implies $\chi_f(G)=3$. One may
ask what is the possible value of $\chi_f(G)$ if $G$ is
triangle-free and  $\Delta(G)$ is at most 3; this problem is
motivated by a well-known and  solved problem of determining the
maximum independence number $\alpha(G)$ for such graphs. Staton
\cite{staton} showed that
\begin{equation}
  \label{eq:2}
\alpha(G) \geq 5n/14
\end{equation}
for any triangle-free graph $G$ on $n$ vertices  with maximum degree
at most $3$. Actually, Staton's bound is the best possible since the
generalized Petersen graph $P(7,2)$ has 14 vertices and independence
number 5 as noticed by Fajtlowicz \cite{fajtlowicz}. Griggs and
Murphy  \cite{gm}  designed a linear-time algorithm to find an
independent set in $G$ of size at least $5(n-k)/14$, where $k$ is
the number of components of $G$ that are 3-regular.  Heckman and
Thomas \cite{ht} gave a simpler proof of Staton's bound and designed
a linear-time algorithm to find an independent set in $G$ with size
at least $5n/14$.

In the same paper \cite{ht}, Heckman and Thomas conjectured
\begin{equation}
  \label{eq:4}
\chi_f(G)\leq \frac{14}{5}
\end{equation}
for every triangle-free graph with maximum degree at most 3.  Note
that \cite{su}
\begin{equation}
  \label{eq:3}
\chi_f(G) = \frac{n}{\alpha(G)},
\end{equation}
provided $G$ is vertex transitive.  Equation (\ref{eq:3}) implies
that the generalized Petersen graph $P(7,2)$ has the fractional
chromatic number $\frac{14}{5}$. Thus, the conjecture is tight if it
holds.

Recently,   Hatami and Zhu  \cite{hz}  proved that $\chi_f(G) \leq 3
-\frac{3}{64}$, provided $G$ is triangle-free with maximum degree at
most three.  Their idea is quite clever. 
For some independent set $X$, the graph obtained by
identifying all neighbors of $X$ into one fat vertex  is $3$-colorable.
Now each vertex in $X$ has freedom of choosing two colors.
This observation results the improvement of $\chi_f(G)$.
Here we improve their result and get
the following theorem.

\begin{theorem}\label{t1}
If $G$ is triangle-free and has maximum degree at most 3, then
$\chi_f(G) \leq 3 -\frac{3}{43}$.
\end{theorem}

\begin{figure}[htb]
    \centering
    \psfig{figure=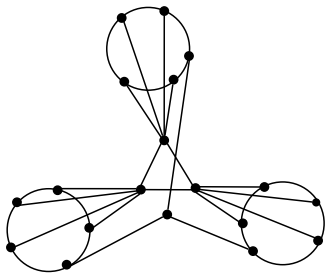, width=0.35\textwidth}
    \label{fig:0}
\caption{A difficult case of $G'(X)$: 
here the triangle in the middle is obtained
from contracting the neighborhood of $X_i$ in $G$, for 
$1\leq i\leq 3$; the remaining
part forms a Gallai tree. Note $\chi(G'(X)) = 4$.}
\end{figure}

Comparing to Hatami and Zhu's result, here we only shrink the gap (to
the conjectured value $2.8$) by $15\%$. However, it is quite hard to
obtain this improvement. The main idea is to extend the independent
set $X$ of $G^\ast$ in Lemma 12 of [5] to the admissible set $X = X_1
\cup X_2\cup X_3$, where $X_1$, $X_2$, and $X_3$ are three independent
sets in $G^\ast$.  However, this extension causes fundamental
difficulty (see the proof of Lemma \ref{l:3colorable}); the reason is
that in general $G'(X)$ (in the proof of Lemma \ref{l:3colorable}) could
be 4-chromatic as shown in Figure \ref{fig:0}.

\begin{figure}[htbp]  
 \centerline{ \psfig{figure=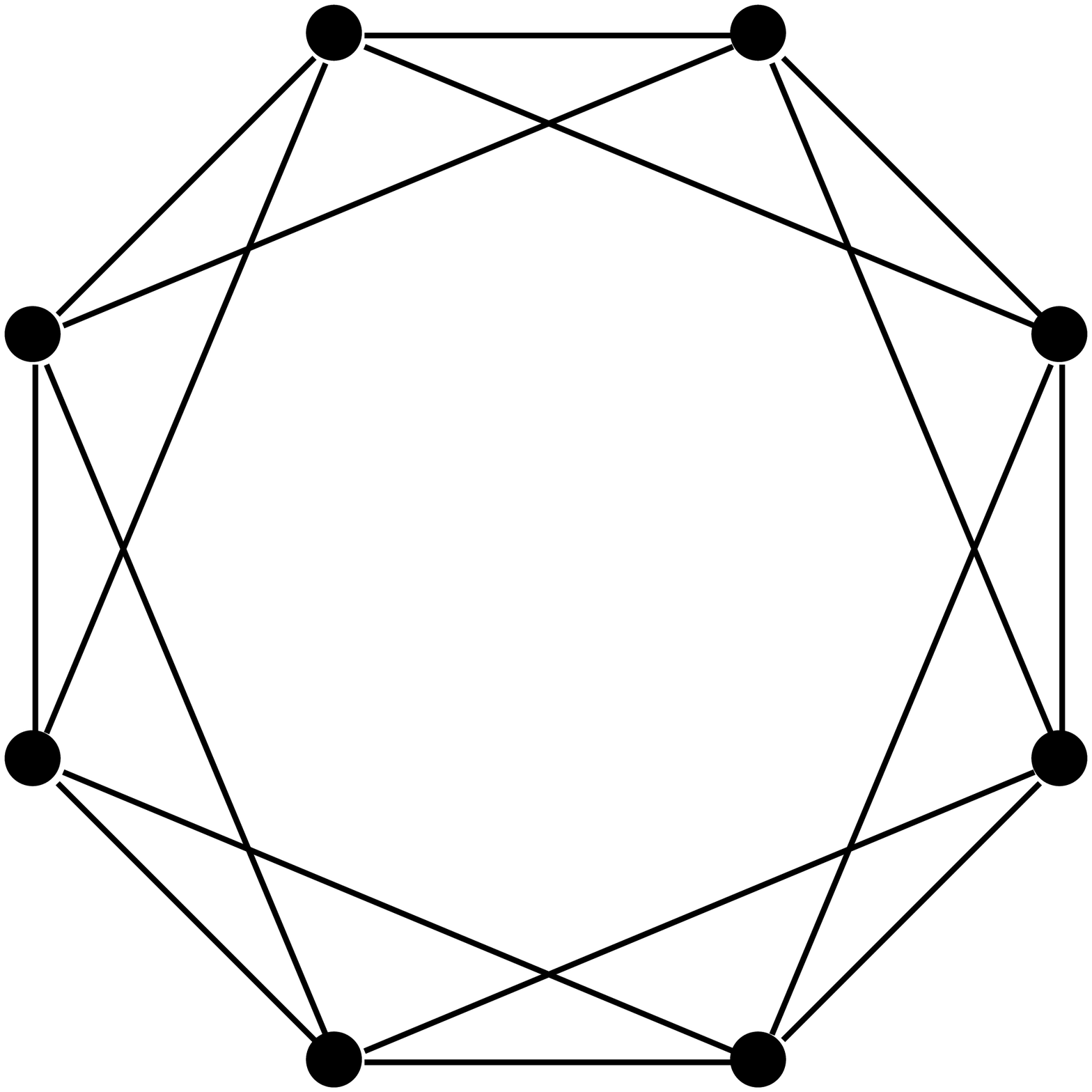, width=0.24\textwidth} 
\hspace{3cm} \psfig{figure=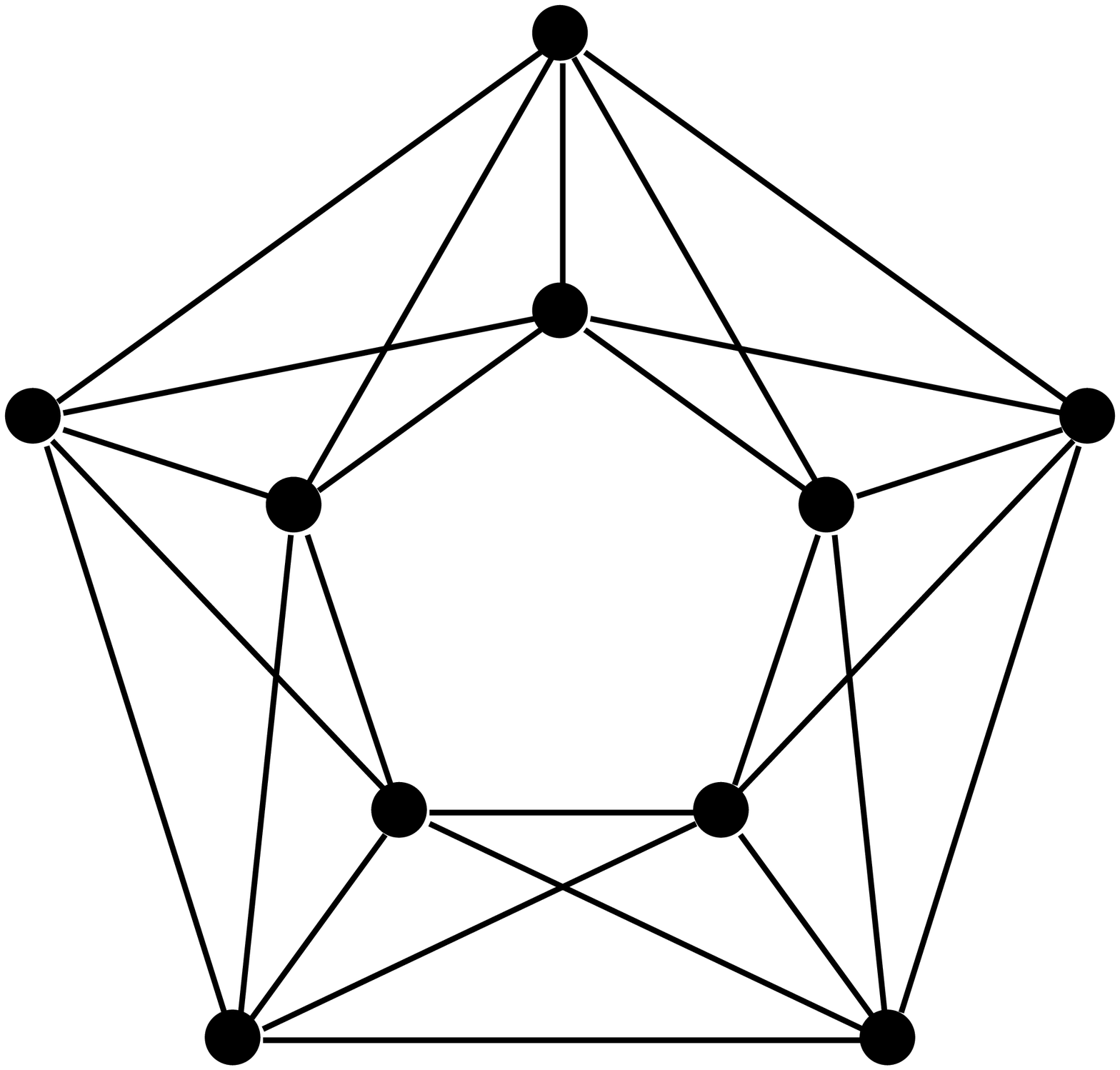,
 width=0.24\textwidth} }
 \centerline{$C_8^2$  \hspace{5cm} $C_5 \boxtimes K_2$}
\caption{Two exception graphs with $\chi_f(G)=\Delta(G)$:
 $C_8^2$ and $C_5 \boxtimes K_2$.}
 \label{fig:2}
\end{figure}

To get over the difficulty, we develop a heavy machinery of
fractionally- critical graphs and use it to prove Lemma \ref{c5} and
\ref{c7}.  The contribution of this paper is not just an improvement
on $\chi_f(G)$ for triangle-free graph $G$ with $\Delta(G)\leq 3$.
More importantly, the theory developed in section 2 can be applied to
more general scenarios concerning the fractional chromatic numbers of
graphs.  Using the tools developed in this paper, King-Lu-Peng
\cite{fbrooks} classified all connected graphs with $\chi_f(G)\geq
\Delta(G)$.  They \cite{fbrooks} further proved that $\chi_f (G) \leq
\Delta(G)- \frac{2}{67}$ for all graphs $G$ such that $G$ is
$K_\Delta$-free, $\Delta(G)\geq 3$, and $G$ is neither $C_8^2$ (the
square of $C_8$) nor $C_5 \boxtimes K_2$ (the strong product of $C_5$
and $K_2$) as shown in Figure \ref{fig:2} Very recently, Edwards and
King \cite{EK} improve the lower bound on $\Delta(G)-\chi_f(G)$ for
all $\Delta\geq 6$.  The Heckman-Thomas' conjecture is a special case
at $\Delta(G)= 3$. We also notice that a better bound $32/11$ (toward
Heckman-Thomas'conjecture) was proved in \cite{FKK} very recently.
 
The rest of the paper is organized as follows. In section 2, we will
study the convex structure of fractional colorings and the
fractionally-critical graphs.  In section 3, we will prove several
key lemmas. In last section, we will show that $G$ can be partitioned
into $42$ admissible sets and present the proof of the main theorem.

\section{Lemmas and Notations}
In this section, we introduce an alternative definition of
``fractional colorings.'' The new definition highlights the convex structure
of the set of all fractional colorings.
The extreme points of these
``fractional colorings'' play a central role in our proofs and seem to have
independent interest.  Our approach is analogous to defining rational
numbers from integers.

\subsection{Convex structures of fractional colorings}
In this paper, we use bold letter $\c$ to represent a coloring.
Recall that a $b$-fold coloring of a graph $G$ assigns a set of $b$
colors to each vertex such that any two adjacent vertices receive
disjoint sets of colors. Given a $b$-fold coloring $\c$, let
$A(\c)=\cup_{v\in V(G)}\c(v)$ be the set of all colors used in $\c$.
Two $b$-fold colorings $\c_1$ and $\c_2$ are {\it isomorphic} if
there is a bijection $\phi\colon A(\c_1) \to A(\c_2)$ such that
$\phi\circ \c_1=\c_2$. In this case we write $\c_1\cong \c_2$. The
isomorphic relation $\cong$ is an equivalence relation. We use $\bar
\c$ to denote the isomorphic class in  which $\c$ belongs to.
Whenever clear under the context, we will not distinguish a $b$-fold
coloring $\c$ and its isomorphic class $\bar \c$.

For a graph $G$ and a positive integer $b$, let $\C_b(G)$ be
the set of all (isomorphic classes of) $b$-fold colorings of $G$.
For  $\c_1\in \C_{b_1}(G)$ and $\c_2\in \C_{b_2}(G)$, we can define
$\c_1+\c_2\in \C_{b_1+b_2}(G)$ as follows: for any $v\in V(G)$,
$$(\c_1+\c_2)(v)=\c_1(v) \sqcup \c_2(v),$$
i.e., $\c_1+\c_2$ assigns $v$ the disjoint union of $\c_1(v)$ and
$\c_2(v)$.

Let $\C(G)=\cup_{b=0}^\infty \C_b(G)$.
It is easy to check that ``$+$'' is commutative and associative.
 Under the addition above, $\C(G)$ forms
a commutative monoid with the unique $0$-fold coloring (denoted by
$0$, for short) as the identity. For a positive integer $t$ and
$\c\in \C_b(G)$, we define
$$t\cdot \c= \overbrace{\c+\cdots+\c}^t$$
to be the new $tb$-fold coloring by duplicating each color $t$ times.

For $\c_1\in \C_{b_1}(G)$ and $\c_2\in \C_{b_2}(G)$, we say $\c_1$ and
$\c_2$ are {\it equivalent}, denoted as $\c_1\sim \c_2$, if there
exists a positive integer $s$ such that $sb_2\cdot \c_1 \cong sb_1
\cdot \c_2$. (This is an analogue of the classical definition of
rational numbers with a slight modification.  
The multiplication by $s$ is needed here 
because the division of a fractional coloring by an integer  has
not been defined yet.
The proofs of Lemma \ref{l:equ} and \ref{l:sim}
are straightforward and are omitted here.)
\begin{lemma}\label{l:equ}
The binary relation $\sim$ is an equivalence relation over $\C(G)$.
\end{lemma}


Let $\F(G)=\C(G)/\sim$ be the set of all equivalence classes. Each
equivalence class is called a fractional coloring of $G$. For any
$\c\in \C_b(G)$, the equivalence class of $\c$ under $\sim$ is
denoted by $\pi(\c)=\frac{\c}{b}$.

\noindent {\bf Remark:} The notation $\frac{\c}{b}$ makes sense only when
$\c$ is a $b$-fold coloring.

Given any rational number $\lambda=\frac{q}{p}\in [0,1]$ (with two positive integers $p$ and $q$) and  two fractional colorings
(two equivalence classes) $\frac{\c_1}{b_1}$ and $\frac{\c_2}{b_2}$,
we define the linear combination as
$$\lambda \c_1 + (1-\lambda) \c_2=\frac{ qb_2 \cdot \c_1 +(p-q)b_1\cdot \c_2}
{pb_1b_2}.$$ The following lemma shows the definition above is
independent of the choices of $\c_1$ and $\c_2$, and so $\lambda
\c_1 + (1-\lambda) \c_2$ is a fractional coloring  depending only on
$\lambda$, $\frac{\c_1}{b_1}$, and $\frac{\c_2}{b_2}$.

\begin{lemma}\label{l:sim}
For $i\in \{1,2,3,4\}$, let  $\c_i \in \C_{b_i}(G)$. Suppose that
$\c_1\sim \c_3$ and $\c_2\sim \c_4$. For any non-negative integers
$p$, $q$, $p'$, and $q'$ satisfying $\frac{q}{p}=\frac{q'}{p'}\in
[0,1]$, we have $qb_2 \cdot \c_1+ (p-q)b_1\cdot \c_2 \sim q'b_4\cdot
\c_3+(p'-q')b_3\cdot \c_4$.
\end{lemma}

Define a function $g_G\colon \F(G)\to {\mathbb Q}$ by
$g_G(\frac{\c}{b})=\frac{|A(\c)|}{b}$. If the graph $G$ is clear
under the context, then we write it as $g(\frac{\c}{b})$ for short.
It is easy to check that $g$ does not depend on the choice of $\c$
and so $g$ is well-defined.  For any $\tau>0$, we define
$$\F_\tau(G)=\left\{\frac{\c}{b}\in \F(G) | g(\frac{\c}{b})\leq \tau\right\}.$$
A fractional coloring $\c$ is called {\em extremal} in $\F_\tau(G)$
if it can not be written as a linear combination of two (or more) distinct 
fractional colorings in  $\F_\tau(G)$.

\begin{theorem} \label{tphi}
For any graph $G$ on $n$ vertices, there is an embedding $\phi\colon
\F(G) \to {\Bbb Q}^{2^n-1}$ such that $\phi$ keeps convex
structure. Moreover, for any rational number $\tau$, $\F_\tau(G)$ is
the convex hull of some extremal fractional colorings.
\end{theorem}
\noindent {\bf Proof:} We would like to define $\phi\colon \F(G)\to {\Bbb Q}^{2^n-1}$
as follows.

Given a fractional coloring $\frac{\c}{b}$,  we can fill these
colors into the regions of the general Venn Diagram on  $n$-sets.
For $1 \leq i \leq  2^n-1$, we can write $i$ as a binary string
$a_1a_2\cdots a_n$ such that  $a_v \in \{0,1\}$ for all $1 \leq
v\leq n$. Write $\c^1(v)=\c(v)$ and $\c^0(v)=\overline{\c(v)}$ (the
complement set of $c(v)$);  then the number of colors in $i$-th
region of the Venn Diagram can be written as
$$h_i(\c)= \left|\cap_{v=1}^n\c^{a_v}(v)
\right|.$$
By the definition, $h_i$ is additive, i.e.
$$h_i(\c_1+\c_2)=h_i(\c_1)+h_i(\c_2).$$
Thus $\frac{h_i(\c)}{b}$  depends only on the fractional coloring $\frac{\c}{b}$
but not on $\c$ itself.

The $i$-th coordinate of $\phi(\frac{\c}{b})$
is defined to be
$$\phi_i(\frac{\c}{b})=\frac{h_i(\c)}{b}.$$
It is easy to check that $\phi_i$ is  a well-defined function on
$\F(G)$.
 Moreover, for any $\lambda=\frac{q}{p}\in [0,1]$
and any two fractional colorings $\frac{\c_1}{b_1}$ and $\frac{\c_2}{b_2}$,
we have
\begin{eqnarray*}
  \phi_i\left(\lambda \frac{\c_1}{b_1} + (1-\lambda) \frac{\c_2}{b_2}\right)
&=&\phi_i\left( \frac{ qb_2 \cdot \c_1 +(p-q)b_1\cdot \c_2}
{pb_1b_2} \right)\\
&=& \frac{h_i(qb_2 \cdot \c_1 +(p-q)b_1\cdot \c_2)}
{pb_1b_2}\\
&=& \frac{qb_2 h_i(\c_1) + (p-q)b_1 h_i(\c_2)}
{pb_1b_2}\\
&=& \frac{q}{p} \frac{h_i(\c_1)}{b_1} + (1-\frac{q}{p})\frac{h_i(\c_2)}{b_2}\\
&=&\lambda \phi_i(\frac{\c_1}{b_1}) + (1-\lambda)\phi_i(\frac{\c_2}{b_2}).
\end{eqnarray*}
Thus $\phi$ keeps the convex structure.

It remains to show $\phi$ is a one-to-one mapping. Assume
$\phi(\frac{\c_1}{b_1})=\phi(\frac{\c_2}{b_2})$. We need to show
$\c_1\sim \c_2$. Let $\c_1'=b_2\cdot \c_1$ and $\c_2'=b_1\cdot
\c_2$. Both $\c_1'$ and $\c_2'$ are $b_1b_2$-fold colorings. Note
$\phi(\frac{\c_1'}{b_1b_2})=\phi(\frac{\c_2'}{b_1b_2})$. For $j \in
\{1,2\}$ and $1\leq i\leq 2^n-1$, we denote the set of colors in the
$i$-th Venn Diagram region of $A(\c_j')$ by $B_i(\c_j')$. Since
$\phi(\frac{\c_1'}{b_1b_2})=\phi(\frac{\c_2'}{b_1b_2})$, we have
$$|B_i(\c_1')| = |B_i(\c_2')|$$
for $1 \leq i \leq 2^n-1.$  There is a bijection $\psi_i$ from
$B_i(\c_1')$ to $B_i(\c_2')$. Note that for $j\in \{1,2\}$,  we have a
partition of $A(\c_j')$:
$$A(\c_j')=\sqcup_{i=1}^{2^n-1} B_i(\c_j').$$
Define a bijection $\psi$
from $A(\c_1')$ to  $A(\c_2')$ be the union of all $\psi_i$
($1\leq i\leq 2^n-1$). We have
$$\psi\circ \c_1'=\c_2'.$$ Thus
$\c_1'\cong \c_2'$. Note that $\c_1\sim \c_1'$ and $\c_2\sim \c_2'$.
We conclude that $\c_1\sim \c_2$.

Under the embedding, $\phi(\F_\tau(G))$ consists of all rational
points in a polytope defined by the intersection of finite number of half
spaces.  Note that all coefficients of the equations of hyperplanes
are rational. Each rational point
in the polytope corresponds to a fractional coloring while each vertex
of the polytope corresponds to an extremal fractional coloring.
\hfill $\square$

\noindent {\bf Remark:} It is well-known that for any graph $G$
there is a $a$:$b$-coloring of $G$ with $\frac{a}{b}=\chi_f(G)$. In
our terminology, we have 
$$\chi_f(G)=\min\left \{g(\frac{\c}{b}) |\
\mbox{for any } \ \frac{\c}{b}\in \F(G) \right \}.$$

\subsection{Coloring restriction and extension}
Let $H$ be a subgraph of $G$. A $b$-fold coloring of $G$ is
naturally a $b$-fold coloring of $H$; this {\it restriction}
operation induces a mapping $i_G^H\colon \F(G)\to \F(H)$. It is easy
to check that $i_G^H$ keeps convex structure,  i.e., for any
$\frac{\c_1}{b_1},\frac{\c_2}{b_2}\in \F(G)$ and $\lambda\in
[0,1]\cap \mathbb{Q}$,   we have
$$i_G^H\left(\lambda \frac{\c_1}{b_1} + (1-\lambda) \frac{\c_2}{b_2}\right)
= \lambda i_G^H\left(\frac{\c_1}{b_1}\right) +
(1-\lambda) i_G^H\left(\frac{\c_2}{b_2}\right).$$
It is also trivial that
$$g_H\left(i_G^H\left(\frac{\c}{b}\right)\right)\leq g_G\left(\frac{\c}{b}\right).$$
Now we consider a reverse operation. We say a fractional coloring
$\frac{\c_1}{b_1}\in \F(H)$ is {\it extensible} in $\F_t(G)$ if
there is a fractional coloring $\frac{\c}{b}\in \F_t(G)$ satisfying
$$ i_G^H\left(\frac{\c}{b}\right)= \frac{\c_1}{b_1}.$$
 We say a fractional coloring
$\frac{\c_1}{b_1}\in \F(H)$ is {\it fully extensible} in $\F(G)$ if
it is extensible in $\F_t(G)$, where  $t= g_H(\frac{\c_1}{b_1})$.
(It also implies that $\frac{\c_1}{b_1}$ is extensible in $\F_t(G)$
for all $t\geq g_H(\frac{\c_1}{b_1})$.)

\begin{lemma}\label{fully}
Let $H$ be a subgraph of $G$ and $ \frac{\c_i}{b_i}\in \F(H)$ for $i
\in \{1,2\}$. Assume that for $i \in \{1,2\}$, $\frac{\c_i}{b_i}$ is
 fully extensible in $\F(G)$.  For any $\lambda\in
{\mathbb Q}\cap [0,1]$,  we have $\lambda \frac{\c_1}{b_1}
+(1-\lambda)\frac{\c_2}{b_2}$  is fully extensible in $\F(G)$.
\end{lemma}
\noindent {\bf Proof:} Let $t_i=g_H(\frac{\c_i}{b_i})$ for $i \in
\{1,2\}$. Note that there are fractional colorings
$\frac{\c_i'}{b_i'} \in \F_{t_i}(G)$ such that $i_G^H \left(
  \frac{\c_i'}{b_i'}\right)=\frac{\c_i}{b_i}$ for $i \in \{1,2\}$.  Let
$\frac{\c}{b}=\lambda\frac{\c_1'}{b_1'} + (1-\lambda)
\frac{\c_2'}{b_2'}$.  We have $g_G(\frac{\c}{b})=\lambda
t_1+(1-\lambda)t_2$ and,
$$
i_G^H\left(\frac{\c}{b} \right) = \lambda
i_G^H\left(\frac{\c_1'}{b_1'}\right) + (1-\lambda)
i_G^H\left(\frac{\c_2'}{b_2'}\right) = \lambda \frac{\c_1}{b_1} +
(1-\lambda) \frac{\c_2}{b_2}.
$$
Note that
$$
g_H\left(\lambda \frac{\c_1}{b_1}+ (1-\lambda)
  \frac{\c_2}{b_2}\right)=\lambda t_1+(1-\lambda)t_2.
$$
Therefore, $\lambda\frac{\c_1}{b_1} +(1-\lambda)\frac{\c_2}{b_2}$ is
fully extensible in $\F(G)$ by the definition. \hfill $\square$

We say $G=G_1\cup G_2$ if $V(G)=V(G_1)\cup V(G_2)$ and $E(G)=E(G_1)
\cup E(G_2)$. Similarly, we say $H=G_1\cap G_2$ if $V(H)=V(G_1)\cap
V(G_2)$ and $E(H)=E(G_1) \cap E(G_2)$.

\begin{lemma}\label{glue}
Let $G$ be a graph.  Assume that $G_1$ and $G_2$ are two subgraphs
such that $G_1 \cup G_2=G$ and $G_1 \cap G_2 =H$.
 If two fractional colorings $\frac{\c_1}{b_1}\in \F(G_1)$
and $\frac{\c_2}{b_2}\in \F(G_2)$ satisfy
$i_{G_1}^H(\frac{\c_1}{b_1})= i_{G_2}^H(\frac{\c_2}{b_2})$, then
there exists a fractional coloring $\frac{\c}{b}\in \F(G)$
satisfying
$$i_G^{G_i}(\frac{\c}{b})= \frac{\c_i}{b_i}$$
for $i \in \{1,2\}$ and
$g_G(\frac{\c}{b})=\max\{g_{G_1}(\frac{\c_1}{b_1}),
g_{G_2}(\frac{\c_2}{b_2})\}$.
\end{lemma}

\noindent {\bf Proof:} Without loss of generality, we can assume
$b_1=b_2=b$ (by taking the least common multiple if it is necessary.) We also assume $g_{G_1}(\frac{\c_1}{b_1})\leq
g_{G_2}(\frac{\c_2}{b_2})$, then we have $|A_{G_1}(\c_1)|\leq
|A_{G_2}(\c_2)|$. Since $i_{G_1}^H(\frac{\c_1}{b_1})=
i_{G_2}^H(\frac{\c_2}{b_2})$, there is a bijection $\phi$ from
$\cup_{v\in V(H)}\c_1(v)$ to $\cup_{v\in
  V(H)}\c_2(v)$.  Extend $\phi$ as an one-to-one mapping from $A_{G_1}(\c_1)$
to $A_{G_2}(\c_2)$ in an arbitrary way. Now we define a $b$-fold
coloring $\c$ of $G$ as follows.
$$\c(v)=\left\{
  \begin{array}{ll}
    \phi(\c_1(v)) & \mbox{ if }v\in V(G_1),\\
    \c_2(v) & \mbox{ if } v\in V(G_2).
  \end{array}
\right.$$ Since $G_1$ and $G_2$ cover all edges of $G$, $\c$ is
well-defined. Note $\c|_{V(G_1)}=\phi \circ \c_1 \cong \c_1$ and
$\c|_{V(G_2)}=\c_2$. Thus for $i \in \{1,2\}$, we have
 $$i_G^{G_i}\left(\frac{\c}{b}\right)= \frac{\c_i}{b_i}.$$
We also have $g_G\left(\frac{\c}{b}\right)=\frac{|A_{G_2}(\c_2)|}{b}
=g_{G_2}(\frac{\c_2}{b_2})$.
\hfill $\square$

\begin{theorem}\label{extension}
Let $G$ be a graph.  Assume that $G_1$ and $G_2$ are two subgraphs
such that $G_1 \cup G_2=G$ and $G_1 \cap G_2 =H$.  If $\chi_f(G_1)
\leq t$ and any extreme fractional coloring in $\F_t(H)$ is
extensible in $\F_t(G_2)$, then we have $\chi_f(G)\leq t$.
\end{theorem}
\noindent {\bf Proof:} There is a fractional coloring
$\frac{\c}{b}\in \F(G_1)$ satisfying $g_{G_1}(\frac{\c_1}{b_1})=t$.
Since every extreme fractional coloring of $\F_t(H)$ is extensible
in $\F_t(G_2)$,  there exist fractional colorings
$\frac{\c_1}{b_1},\frac{\c_2}{b_2}, \cdots, \frac{\c_r}{b_r} \in
\F_t(G_2)$ such that
$i_{G_2}^H(\frac{\c_1}{b_1}),i_{G_2}^H(\frac{\c_2}{b_2}),\cdots,
i_{G_2}^H(\frac{\c_r}{b_r})$ are all extreme fractional colorings in
$\F_t(H)$. The fractional coloring $i_{G_1}^H(\frac{\c}{b})$ can be
written as a linear combination of the extreme ones. Hence, there
exist $\lambda_1,  \lambda_2,\ldots, \lambda_r \in \mathbb{Q}\cap
[0,1]$ such that $\sum_{i=1}^r \lambda_i=1$ and
$$i_{G_1}^H(\frac{\c}{b})=\sum_{i=1}^r \lambda_i i_{G_2}^H(\frac{\c_i}{b_i}).$$
Let $\frac{\c'}{b'}=\sum_{i=1}^r  \lambda_i\frac{\c_i}{b_i} \in \F_t(G_2)$.
We have
$$i_{G_1}^H(\frac{\c}{b})=i_{G_2}^H(\frac{\c'}{b'}).$$
Applying Lemma \ref{glue}, there exists a fractional coloring
$\frac{\c''}{b''}\in \F_t(G)$. Therefore, we have $\chi_f(G)\leq t$.
\hfill $\square$

\begin{corollary}\label{cor1}
Let $G$ be a graph.  Assume that $G_1$ and $G_2$ are two subgraphs
such that $G_1 \cup G_2=G$ and $G_1 \cap G_2 =K_r$ for some positive
integer $r$. We have
$$\chi_f(G)=\max\{\chi_f(G_1), \chi_f(G_2)\}.$$
\end{corollary}
\noindent {\bf Proof:} Without loss of generality, we assume
$\chi_f(G_1)\geq \chi_f(G_2)$. Let $t=\chi_f(G_1)$. We have $t\geq
r$ as $G_1$ contains $K_r$. Since $\chi_f(G_2)\leq t$,  we have
$\F_t(G_2)\not=\emptyset$. Since $H$ is a complete graph, $\F_t(H)$
contains only one fractional coloring; namely, color all vertices of
$H$ using distinct colors. It is trivial that any extreme fractional
coloring in $\F_t(H)$ is extensible in $\F_t((G_2)$. Applying Theorem
\ref{extension}, we have $\chi_f(G)\leq t$. The other direction is
trivial. \hfill $\square$

Let $uv$ be a non-edge of a graph $G_2$. We denote $G_2+uv$ to be
the supergraph of $G_2$ by adding the edge $uv$ and denote $G_2/uv$
be the quotient graph by identifying the vertex $u$ and the vertex
$v$.

\begin{figure}[htbp]
 \centerline{ {\psfig{figure=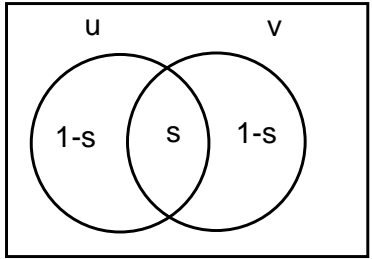, width=0.40\textwidth}}}
\caption{Fractional colorings on the empty graph on two vertices $u$ and $v$
are represented as a weighted Venn Diagram.} \label{ven}
\end{figure}

\begin{lemma}\label{l:cut2}
Let $G$ be a graph.  Assume that $G_1$ and $G_2$ are two subgraphs
such that $G_1 \cup G_2=G$ and $V(G_1) \cap V(G_2) =\{u,v\}$.
\begin{enumerate}
\item If $uv$ is an edge of $G$, then we have
$$\chi_f(G) = \max\{\chi_f(G_1), \chi_f(G_2)\}.$$
\item If $uv$ is not an edge of $G$, then we have
$$\chi_f(G) \leq \max\{\chi_f(G_1), \chi_f(G_2+uv), \chi_f(G_2/uv)\}.$$
\end{enumerate}
\end{lemma} \noindent {\bf Proof:} Part 1 is a simple application of Corollary
\ref{cor1}.  For the proof of part 2, let $t= \max\{\chi_f(G_1),
\chi_f(G_2+uv), \chi_f(G_2/uv)\}$. Note $t\geq 2$. 
Let $E_{uv}$ be the empty graph on the set of two vertices $u$ and $v$.
All fractional
colorings of $\F_t(E_{uv})$ can be represented by the following
weighted Venn Diagram, see Figure \ref{ven}.

The parameter $s$ measure the fraction of common colors shared by
$u$ and $v$. There are two extreme points in the convex hull: $s=0$ and $s=1$.  
The extremal fractional coloring corresponding to $s=0$ is extensible in
$\F_t(G_2)$ since $\chi_f(G+uv) \leq t$. The extremal fractional coloring
corresponding to $s=1$ is extensible in $\F_t(G_2)$ since
$\chi_f(G/uv) \leq t$.  Applying Theorem \ref{extension} and Lemma 
\ref{fully}, we have
$\chi_f(G)\leq t$. Part 2 is proved. \hfill $\square$

\subsection{Fractionally-critical graphs}
In this subsection, we will
apply our machinery to triangle-free graphs with maximum degree at
most $3$.

Recall that a graph $G$ is {\it $k$-critical} (for a positive
integer $k$) if $\chi(G)=k$ and $\chi(H)<k$ for any proper subgraph
$H$ of $G$.  For any rational number $t\geq 2$, a graph $G$ is {\it
$t$-fractionally-critical} if $\chi_f(G)=t$ and $\chi_f(H)<t$ for
any proper subgraph $H$ of $G$. For simplicity, we say $G$ is {\it
fractionally-critical} if $G$ is $\chi_f(G)$-fractionally-critical.

We will study the properties of  fractionally-critical
graphs. The following lemma is a consequence of Corollary \ref{cor1}.

\begin{lemma}
Assume that $G$ is  a fractionally-critical graph with
$\chi_f(G)\geq 2$. We have $G$ is $2$-connected. Moreover, if   $G$
has a vertex-cut $\{u,v\}$, then $uv$ is not an edge of $G$.
\end{lemma}

For any vertex $u$ of a graph $G$ and a positive integer $i$, we define\\
\centerline{$N_G^i(u)=\{v\in V\colon v\not=u \mbox{ and there is a path of
length $i$ connecting } u \mbox{ and }v \}$.}

\begin{lemma}\label{c5}
Assume that $G$ is a fractionally-critical triangle-free graph satisfying  $\Delta(G) \leq 3$ and $\frac{11}{4} < \chi_f(G) < 3$. For any
vertex $x\in V(G)$ and any $5$-cycle $C$ of $G$, we have either $|V(C)\cap
N_G^2(x)|\leq 3$ or $|V(C)\cap N_G^1(x)|\geq 1$.
\end{lemma}
\noindent {\bf Proof:} Let $t=\chi_f(G)$. We have $\frac{11}{4}<t <
3$. We will prove the statement by contradiction. Suppose that there
is a vertex $x$ and a $5$-cycle $C$ satisfying $|V(C)\cap
N_G^2(x)|\geq 4$ and $V(C)\cap  N_G^1(x)=\emptyset$. Combined with the fact $G$ being triangle-free, we
have the following two cases.

\begin{figure}[htbp]
  \centerline{ \psfig{figure=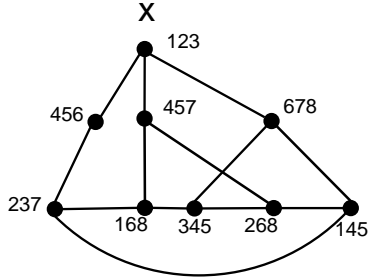, width=0.4\textwidth}}
\caption{An $8\!:\!3$-coloring of $G_9$, where $|V(C)\cap
N_G^2(x)|=5$.} \label{fig:c5g9}
\end{figure}

\noindent
{\bf Case 1:}
$|V(C)\cap N_G^2(x)|=5$. Since $\Delta(G)\leq 3$ and $G$ is triangle-free,
it is easy to check that $G$ contains
  the following subgraph $G_9$ as shown in Figure \ref{fig:c5g9}.
  Since $G$ is $2$-connected, $G_9$ is the entire graph (in
  \cite{hz}). Thus $\chi_f(G)\leq \frac{8}{3}<t$. Contradiction!

\begin{figure}[htbp]
 \centerline{ {\psfig{figure=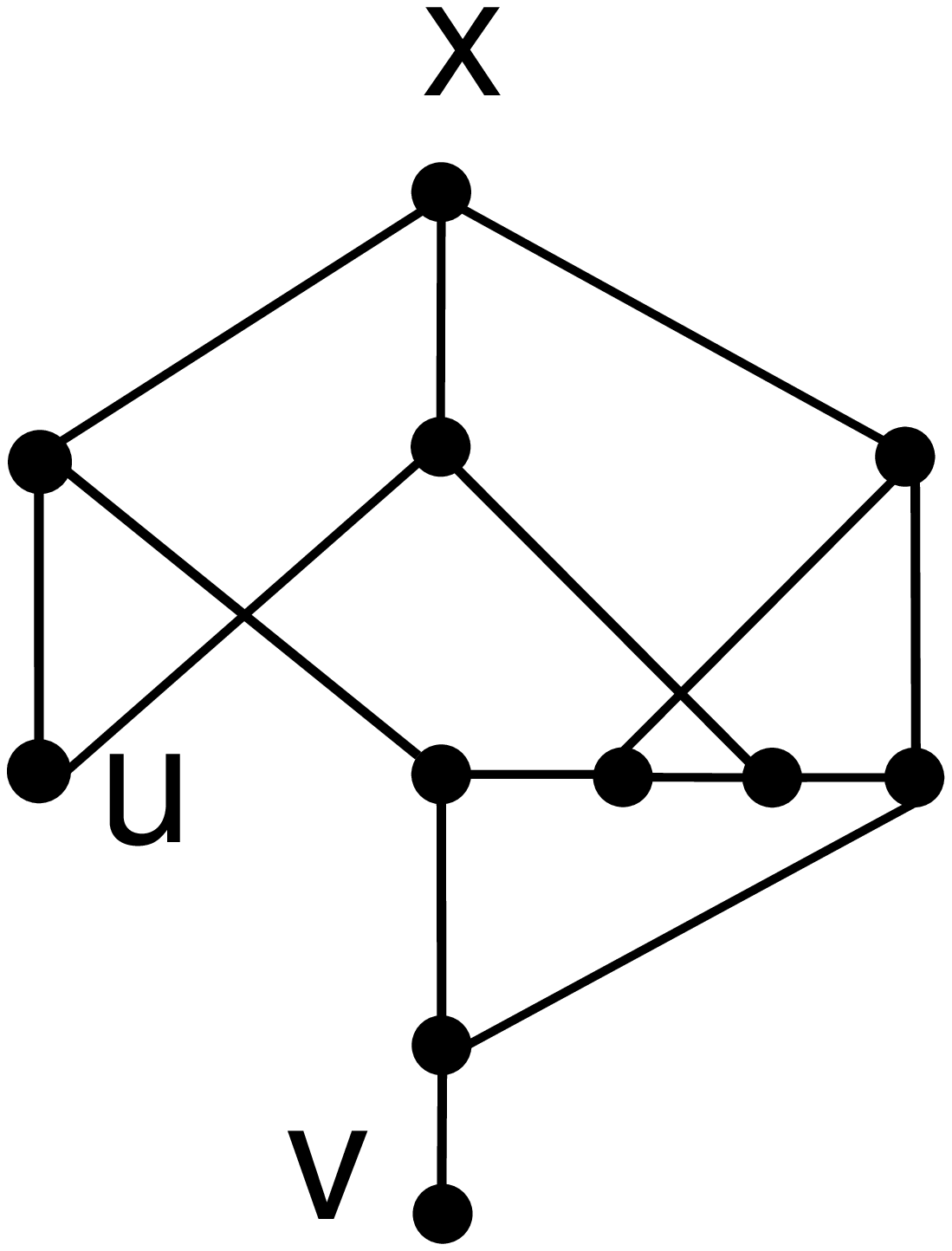, width=0.33\textwidth}}
    \hfil \psfig{figure=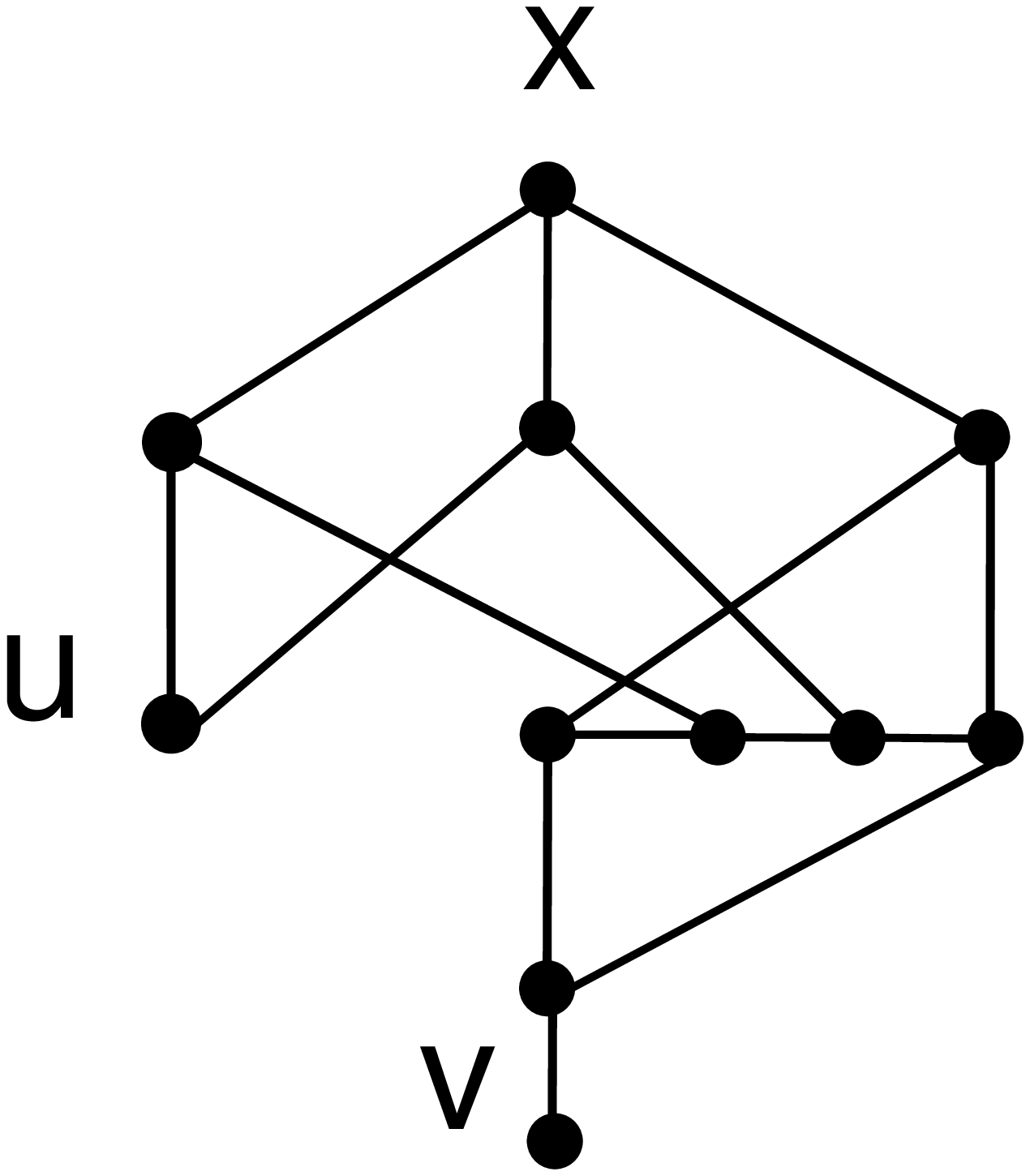, width=0.33\textwidth}
\hfil \psfig{figure=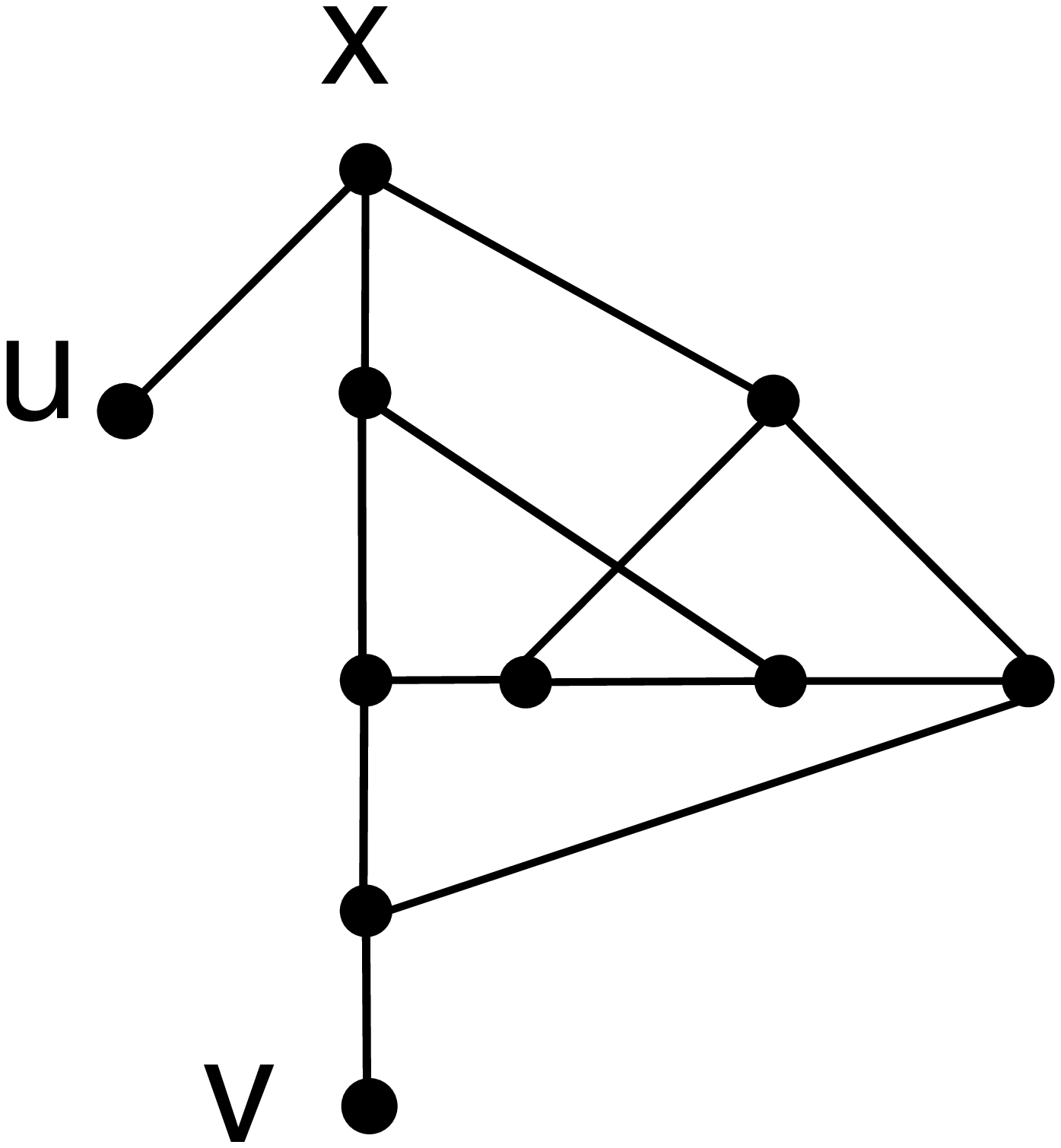, width=0.33\textwidth}}
\centerline{ (I)\hspace*{1cm}\hfil\hspace*{1cm} (II)\hspace*{1cm} \hfil\hspace*{1cm} (III)}
  \centerline{ {\psfig{figure=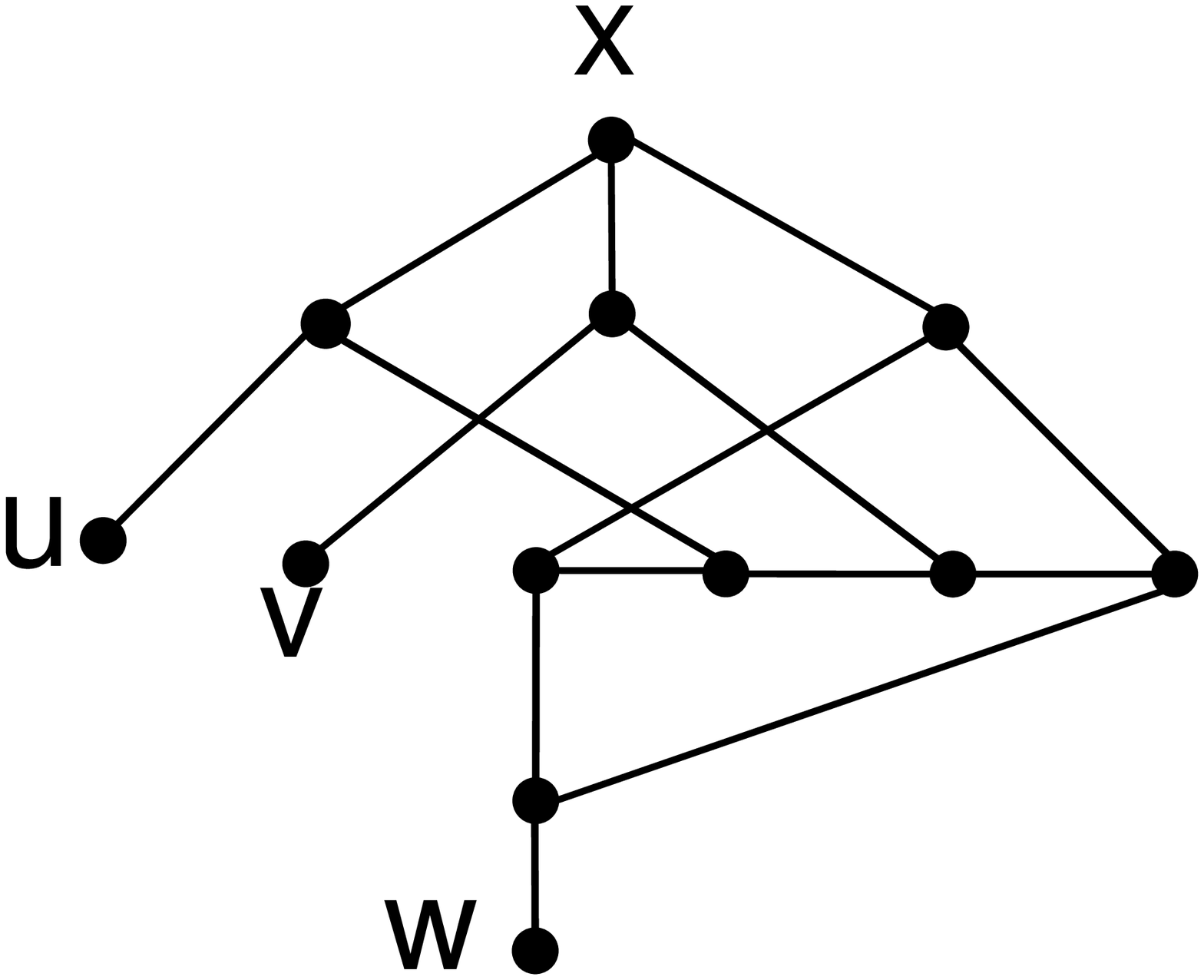, width=0.33\textwidth}}
    \hfil \psfig{figure=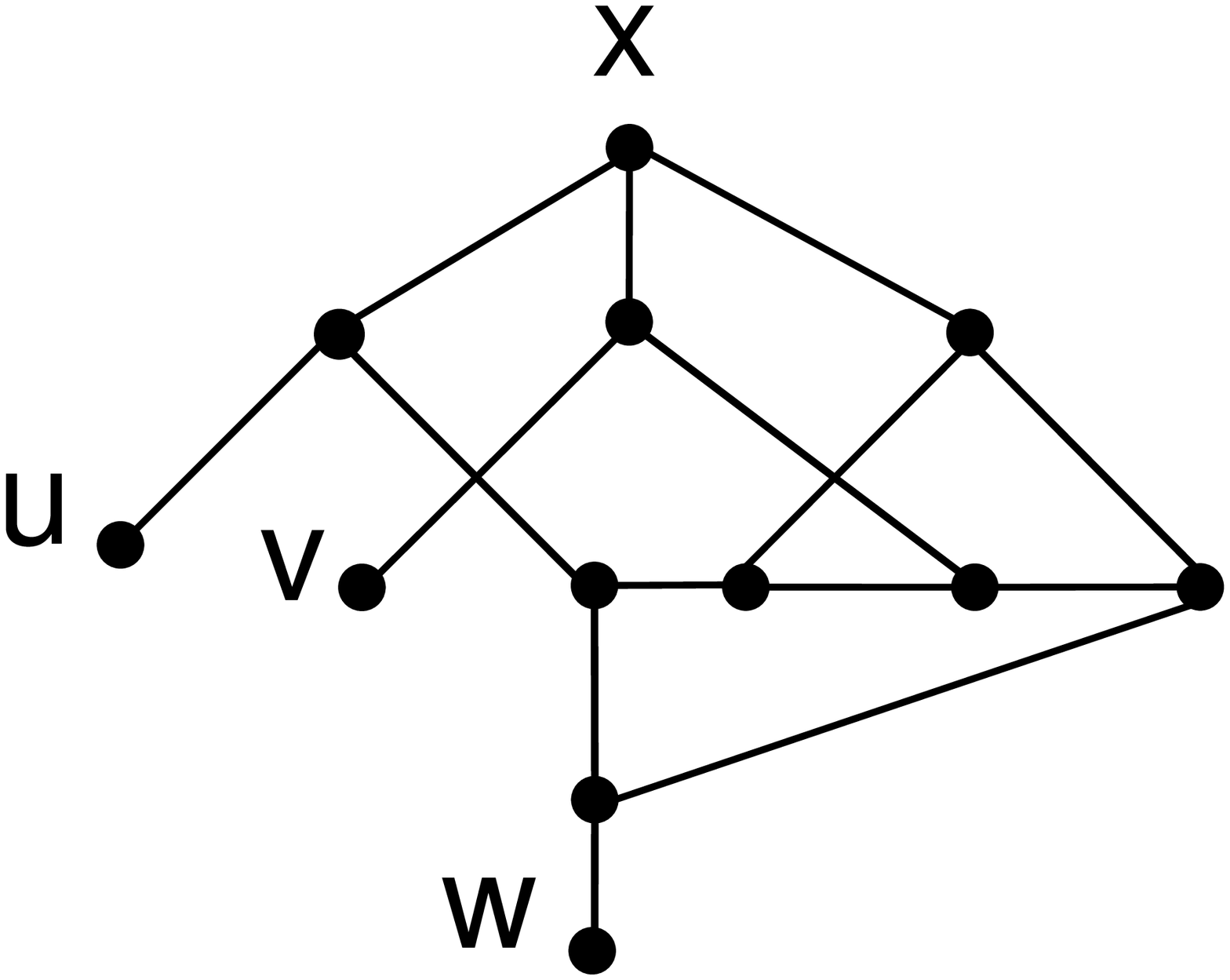, width=0.33\textwidth}}
\centerline{ (IV) \hspace*{1cm}\hfil \hspace*{1cm}(V)}
\caption{All possible cases of $|V(C)\cap N_G^2(x)|=4$ and
$|V(C)\setminus (N_G^1(x)\cup N_G^2(x))|=1$.}
\label{fig:c5s}
\end{figure}

\noindent {\bf Case 2:} $|V(C)\cap N_G^2(x)|=4$ and there exists one
vertex of $C$ having  distance of $3$ to $x$. Hatami and Zhu
\cite{hz}  showed that
 $G$ contains one of the five graphs in Figure \ref{fig:c5s} as a subgraph.
(Note that some of the marked vertices $u$, $v$, and $w$ may be
missing or overlapped;  these degenerated cases result in a smaller
vertex-cut, and can be covered in a similar but easier way; we will
discuss them at the end of this proof.)

\begin{figure}[htbp]
  \centerline{ {\psfig{figure=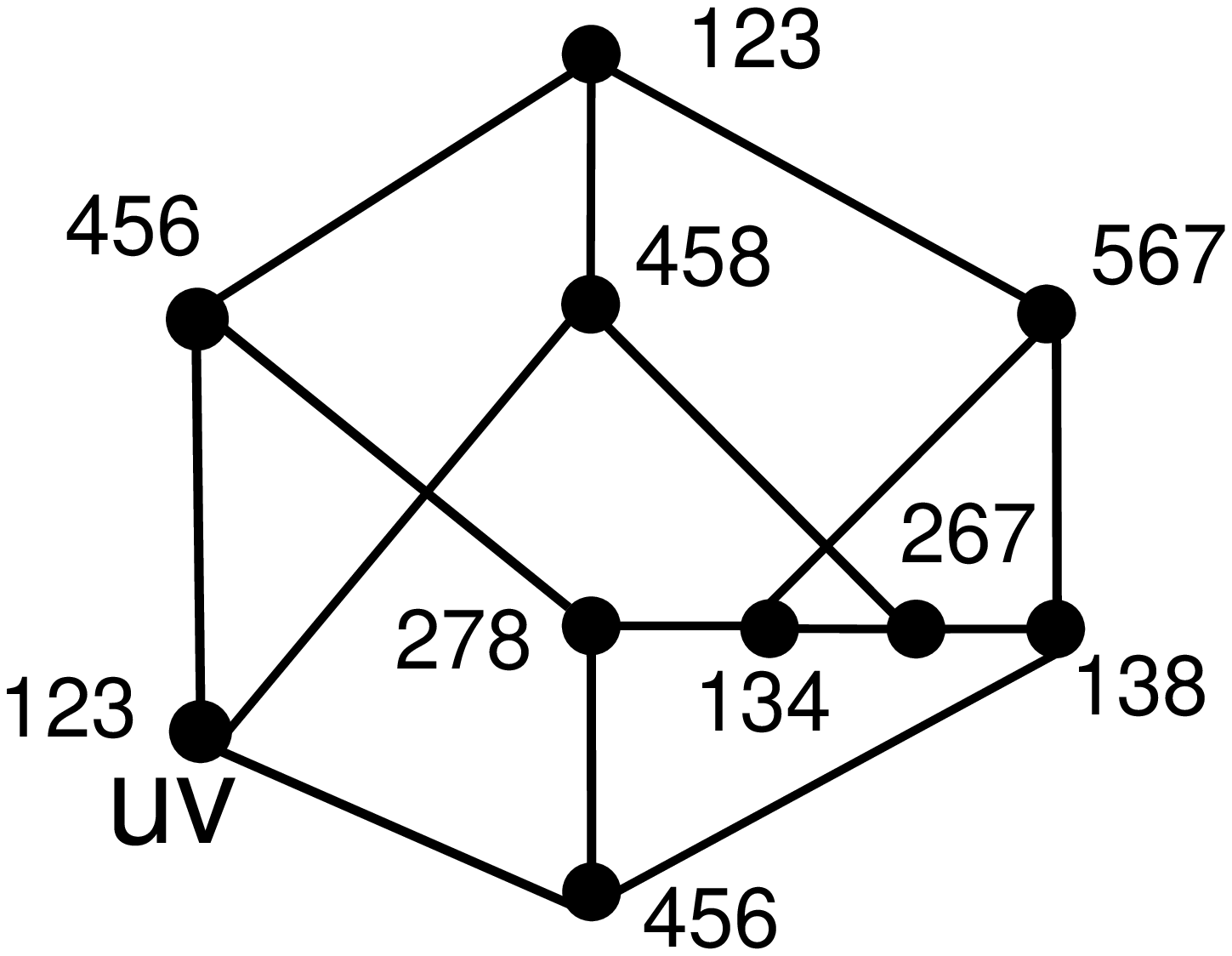,
        width=0.33\textwidth}} \hfil \psfig{figure=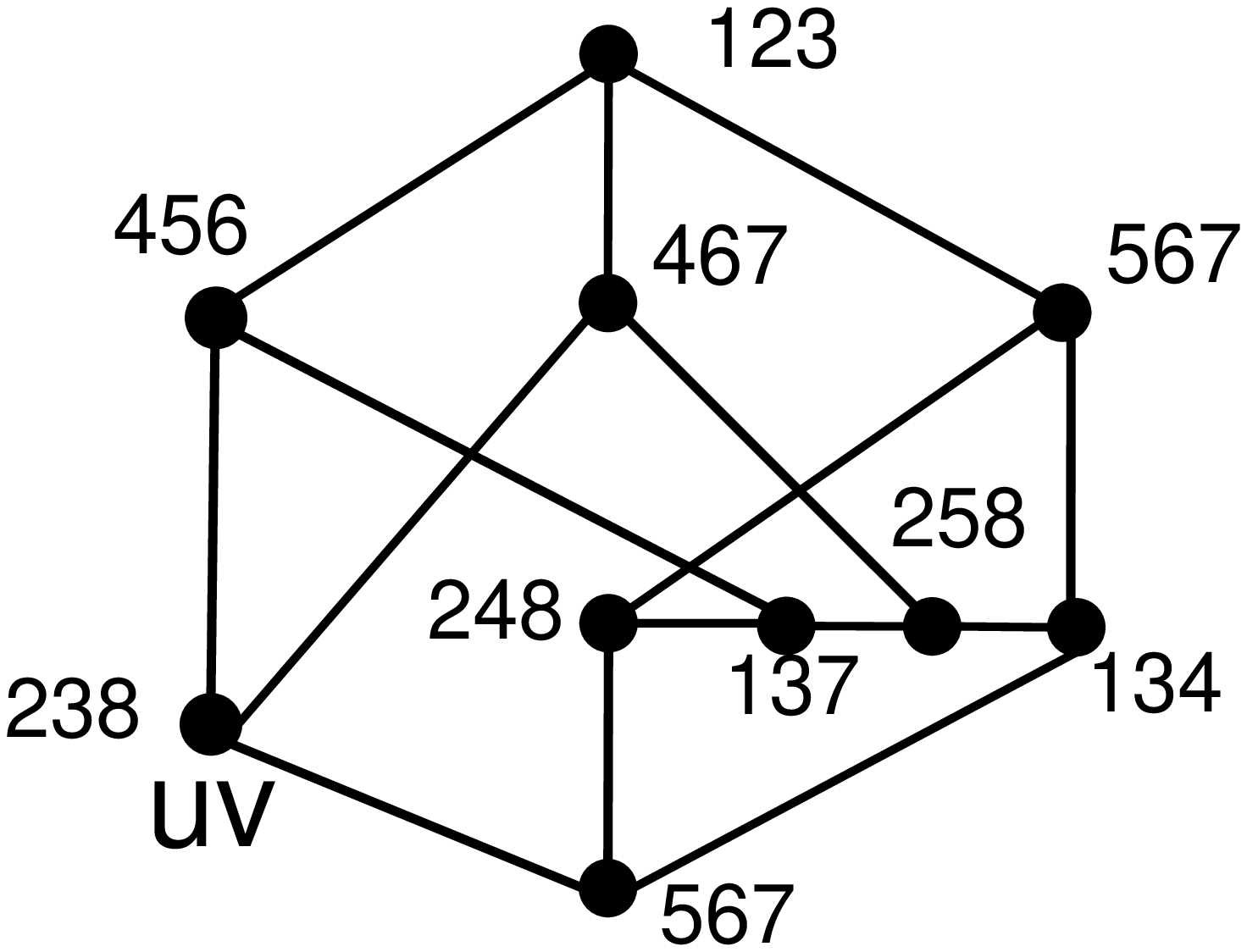,
      width=0.33\textwidth} \hfil \psfig{figure=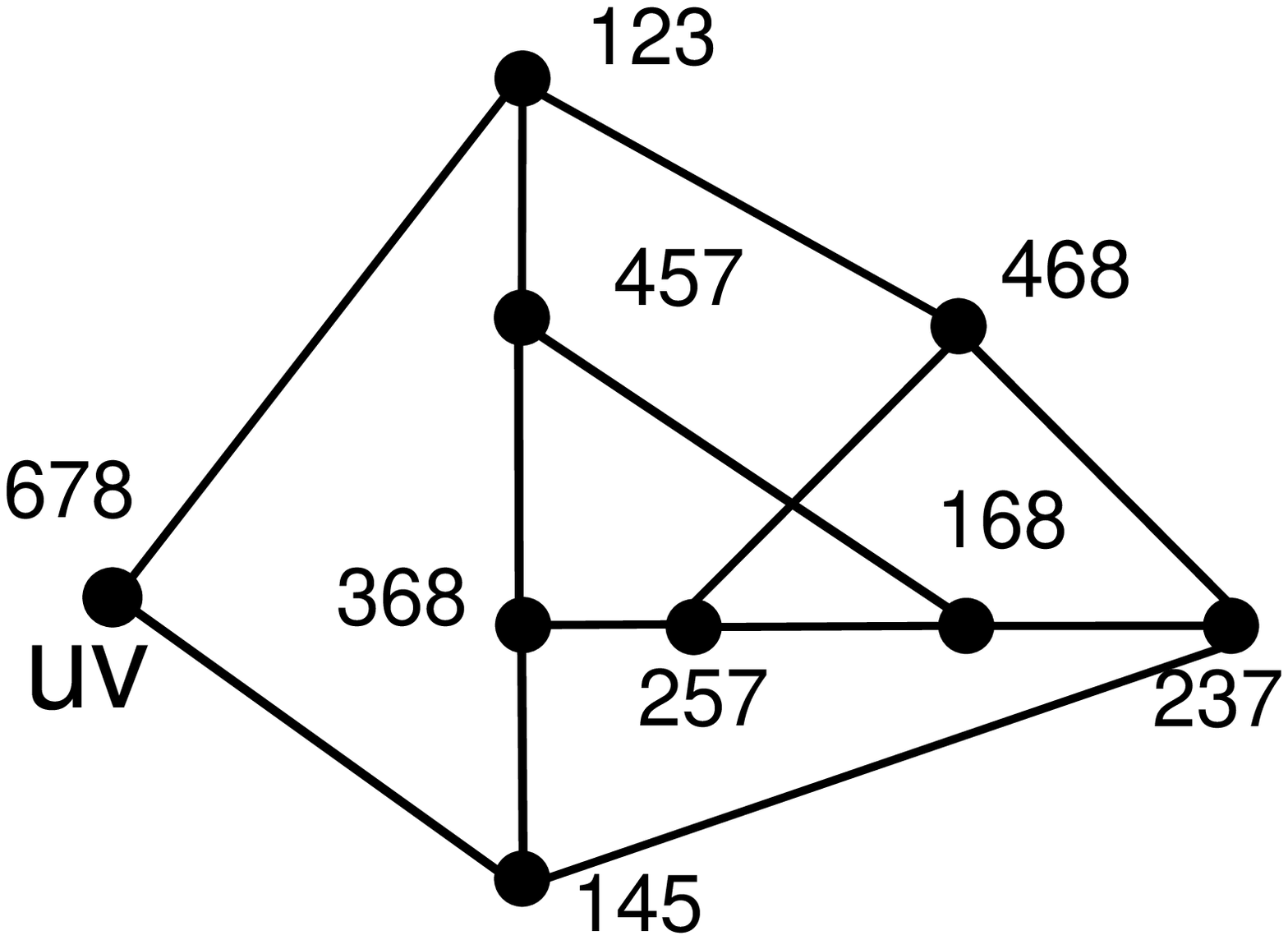,
      width=0.33\textwidth}} \centerline{ $G_2/uv$ in (I)
    \hspace*{1cm}\hfil\hspace*{1cm}$G_2/uv$ in (II)\hspace*{1cm}
    \hfil\hspace*{1cm} $G_2/uv$ in (III)} \centerline{
    {\psfig{figure=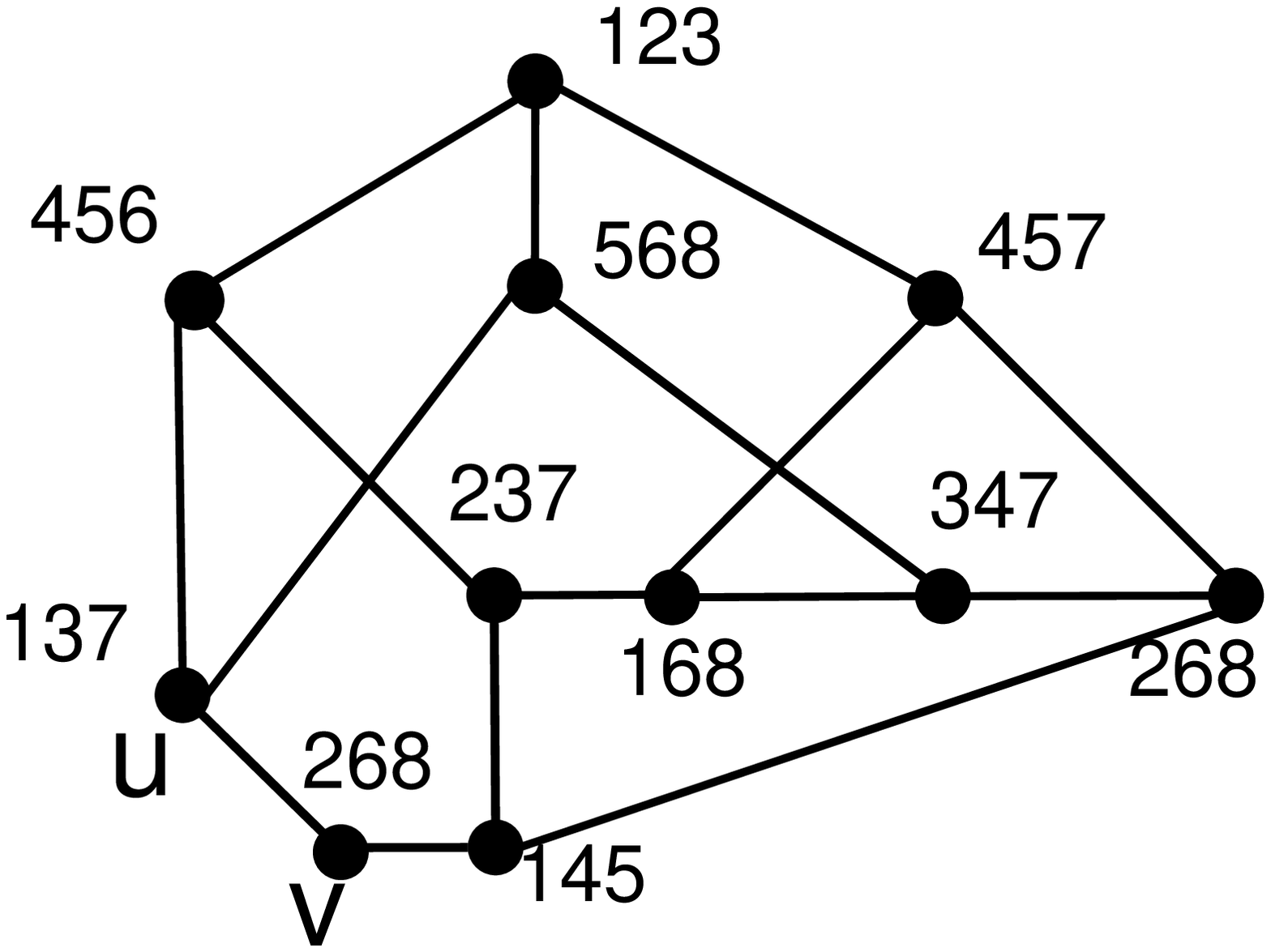, width=0.33\textwidth}} \hfil
    \psfig{figure=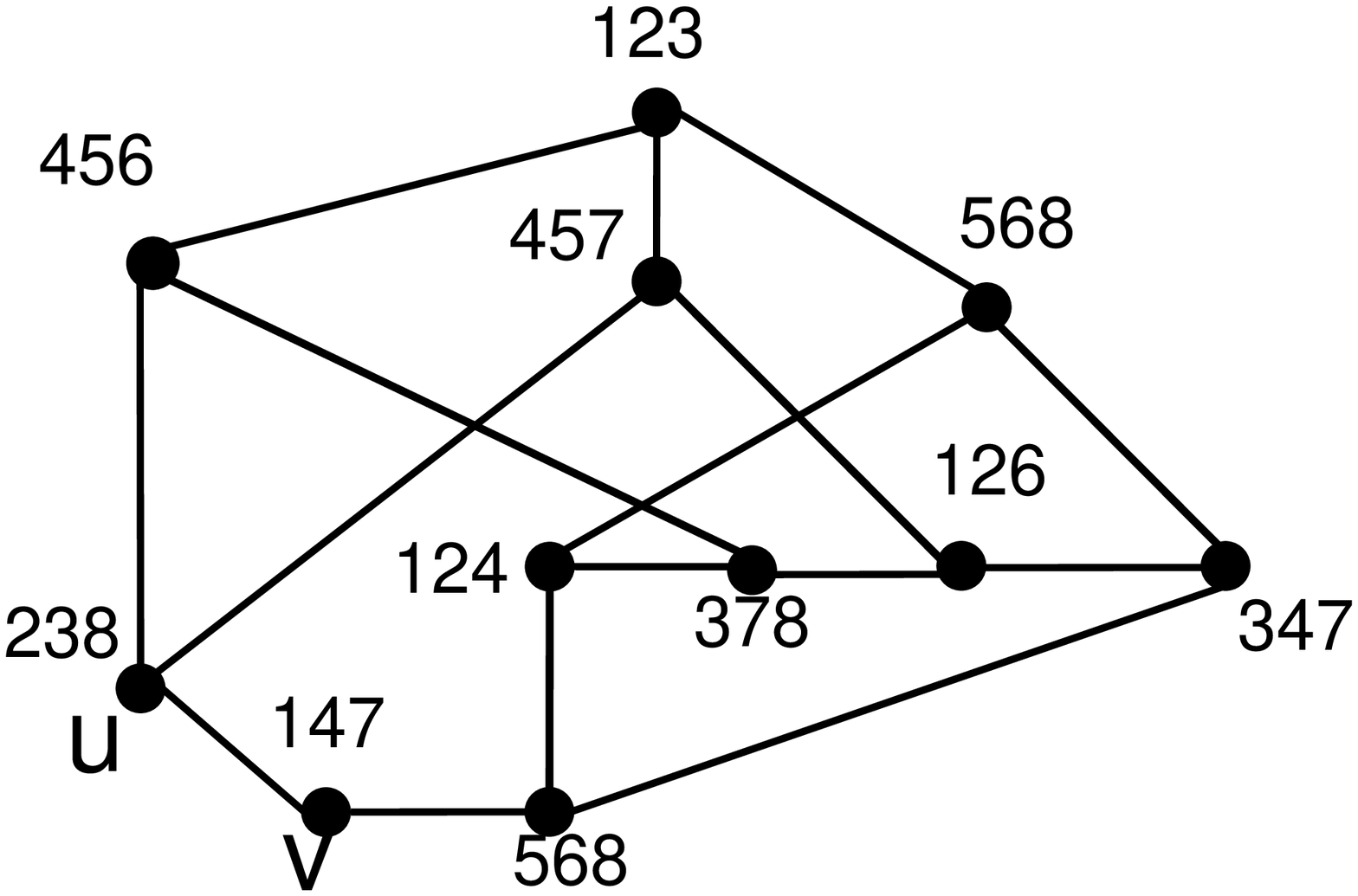, width=0.33\textwidth} \hfil
    \psfig{figure=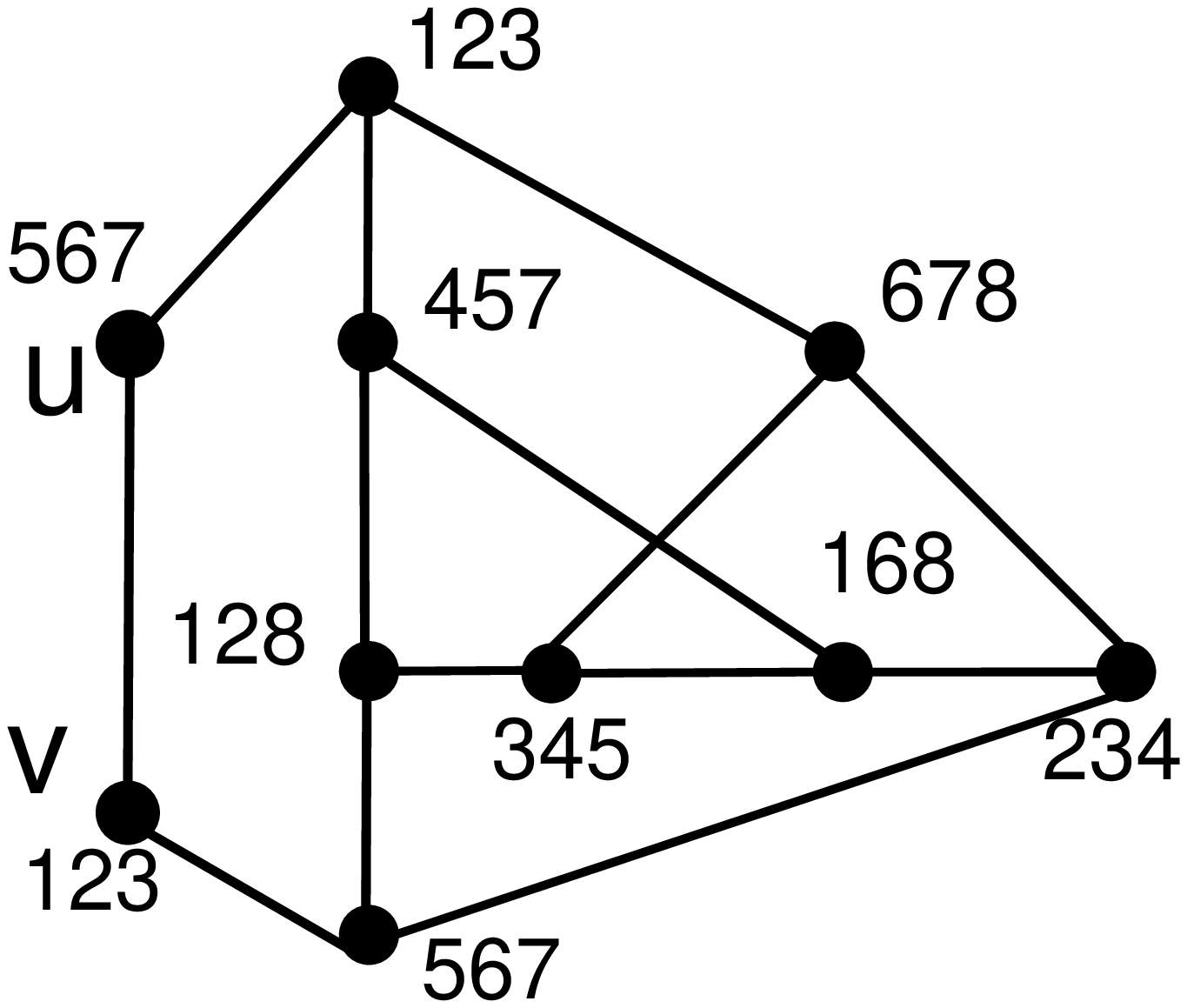, width=0.33\textwidth}}
  \centerline{ $G_2+uv$ in (I) \hspace*{1cm}\hfil\hspace*{1cm}$G_2+uv$
    in (II)\hspace*{1cm} \hfil\hspace*{1cm} $G_2+uv$ in (III)}
\caption{$G_2+uv$ and $G_2/uv$ are all $8\!:\!3$-colorable in cases (I), (II), and (III).
}
\label{fig:c5s2}
\end{figure}

If $G$ contains a subgraph of type (I), (II), or (III), then $G$ has
a  vertex-cut $\{u, v\}$.  Let $G_1$ and $G_2$ be the two connected
subgraphs of $G$ such that $G_1\cup G_2=G$, $V(G_1)\cap
V(G_2)=\{u,v\}$, and $x\in G_2$. In all three cases, $G_2+uv$ and
$G_2/uv$ are $8\!:\!3$-colorable. Please see Figure \ref{fig:c5s2}.

Applying Lemma \ref{l:cut2}, we have
\begin{eqnarray*}
\chi_f(G)&\leq& \max\{\chi_f(G_1), \chi_f(G_2/uv), \chi_f(G_2+uv)\}\\
&\leq& \max\left\{\chi_f(G_1), \frac{8}{3}\right\}.
\end{eqnarray*}
Since $\chi_f(G)\geq t>\frac{8}{3}$,  we must have $\chi_f(G_1)\geq
\chi_f(G)=t$, which is a contradiction to the assumption that $G$
is fractionally-critical.

\begin{figure}[htbp]
  \centerline{\psfig{figure=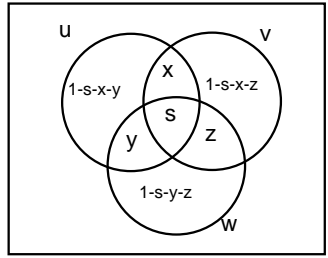, width=0.35\textwidth}}
 \caption{The general fractional colorings on the vertices $u$, $v$, and $w$.}
\label{fig:venn1}
\end{figure}

If $G$ contains one of the subgraphs (IV) and (V), then $G$ has a
vertex-cut set $H=\{u,v,w\}$ as shown in Figure \ref{fig:c5s}. Let
$G_1$ and $G_2$ be the two connected subgraphs of $G$ such that
$G_1\cup G_2=G$, $G_1\cap G_2=\{u,v,w\}$, and $x\in G_2$. We shall
show $\chi_f(G_1)=t=\chi_f(G)$.  Suppose not. We can assume
$\chi_f(G_1) \leq t_0 < t$, where $\frac{8}{3} < \frac{11}{4} <
t_0<t<3$.  By Theorem \ref{tphi}, every fractional coloring in $\F_{t_0}(H)$ 
can be represented by a rational point in a convex polytope $\phi(\F_{t_0}(H))$.
From now on, we will not distinguish the rational point in the convex polytope
and the fractional coloring. Note the convex polytope for
$\F_{t_0}(H)$  can be parametrized as
$$\F_{t_0}(H) \cong \left\{(x,y,z,s)\in {\mathbb Q}^4\left |
  \begin{array}{c}
x+y+s\leq 1\\
x+z+s\leq 1\\
y+z+s\leq 1\\
3-x-y-z-2s\leq t_0 \\
x,y,z,s\geq 0
  \end{array}
\right. \right\}.$$ See the weighted Venn Diagram in Figure
\ref{fig:venn1}.

The extreme fractional colorings of $\F_{t_0}(H)$ are represented (under $\phi$, see Theorem \ref{tphi}) by:
\begin{description}
\item[(a)]  $x=y=z=0$ and $s=1$.
\item[(b)]  $x=1$ and $y=z=s=0$.
\item[(c)]  $y=1$ and $x=z=s=0$.
\item[(d)]  $z=1$ and $x=y=s=0$.
\item[(e)]  $x=3-t_0$ and $y=z=s=0$.
\item[(f)]  $y=3-t_0$ and $x=z=s=0$.
\item[(g)]  $z=3-t_0$ and $x=y=s=0$.
\end{description}
We will show that all 7 extreme fractional colorings
are extensible in $\F_{t_0}(G_2)$.
\begin{description}
\item[(a)] Let $G_2/uvw$ be the quotient graph
by identifying $u$, $v$, and $w$ as one vertex. The fractional
coloring $(0,0,0,1)$ is extensible in $\F_{t_0}(G_2)$ if and only if
$\chi_f(G_2/uvw)\leq t_0$, which is verified by Figure
\ref{fig:uvw}.
\begin{figure}[htbp]
  \centerline{ {\psfig{figure=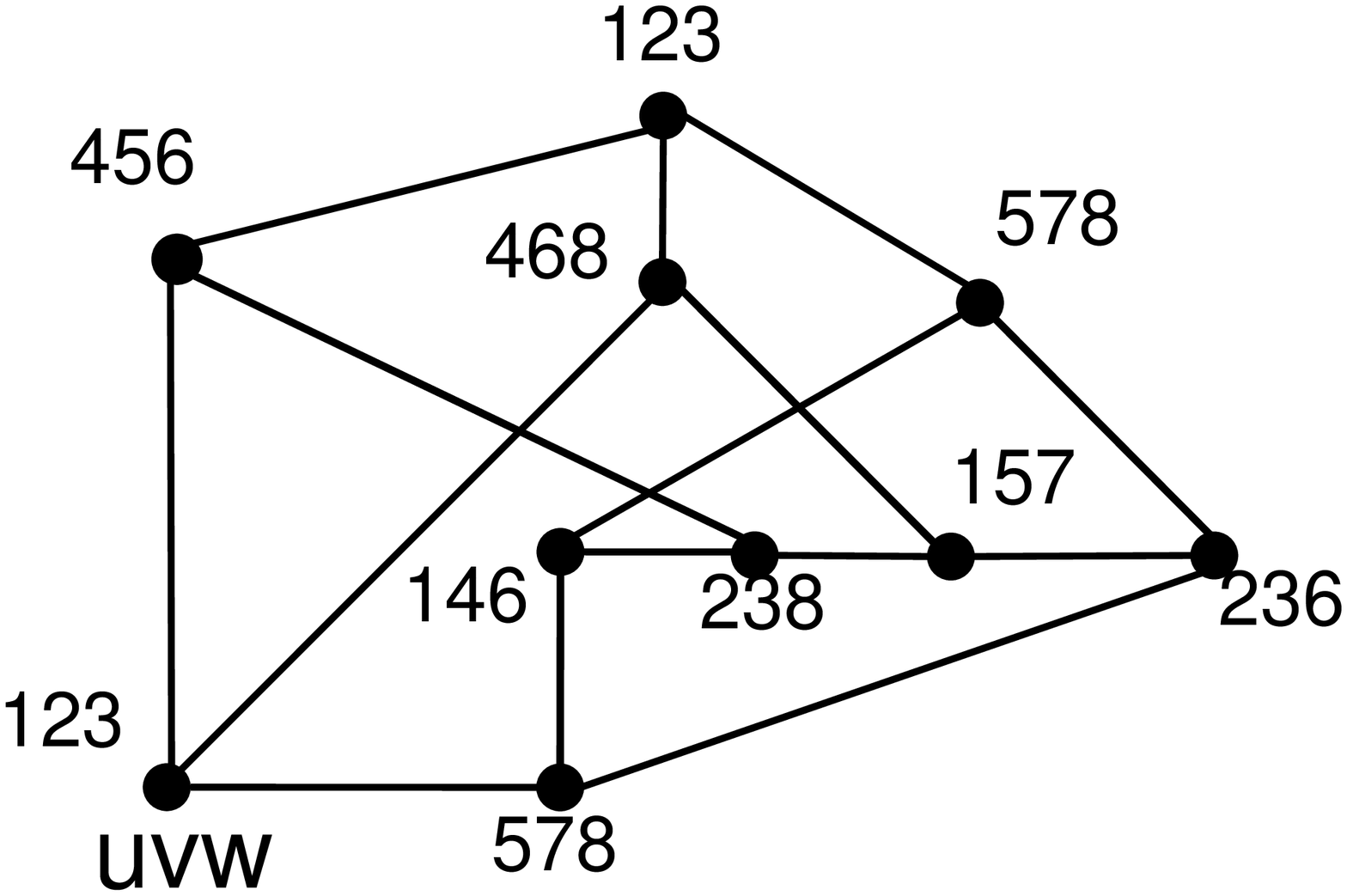, width=0.4\textwidth}}
    \hfil \psfig{figure=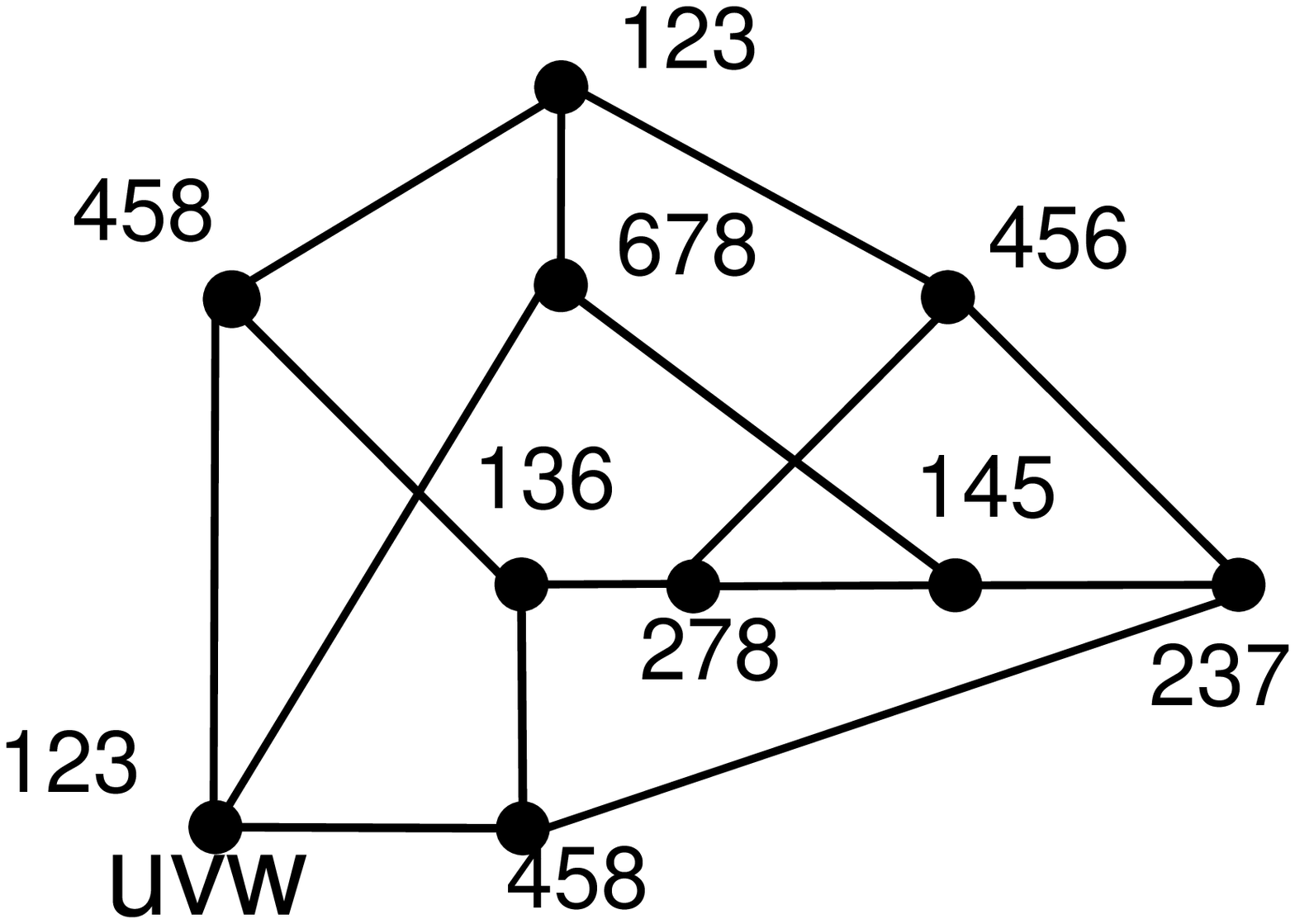, width=0.4\textwidth}}
\centerline{$G_2/uvw$ in (IV) \hspace*{1cm} \hfil \hspace*{1cm}$G_2/uvw$ in (V)}
\caption{$G_2/uvw$ are  $8\!:\!3$-colorable for subgraphs (IV) and (V).}
\label{fig:uvw}
\end{figure}
\item[(b)] Let $(G_2/vw)+ u(vw)$ be the graph obtained by
identifying $v$  and $w$ as one vertex $vw$ followed by adding an
edge $u(vw)$. The fractional coloring $(0,0,1,0)$ is extensible in
$\F_{t_0}(G_2)$ if and only if $\chi_f((G_2/vw)+ u(vw))\leq t_0$,
which is verified by Figure \ref{fig:u-vw}.
\begin{figure}[htbp]
  \centerline{ {\psfig{figure=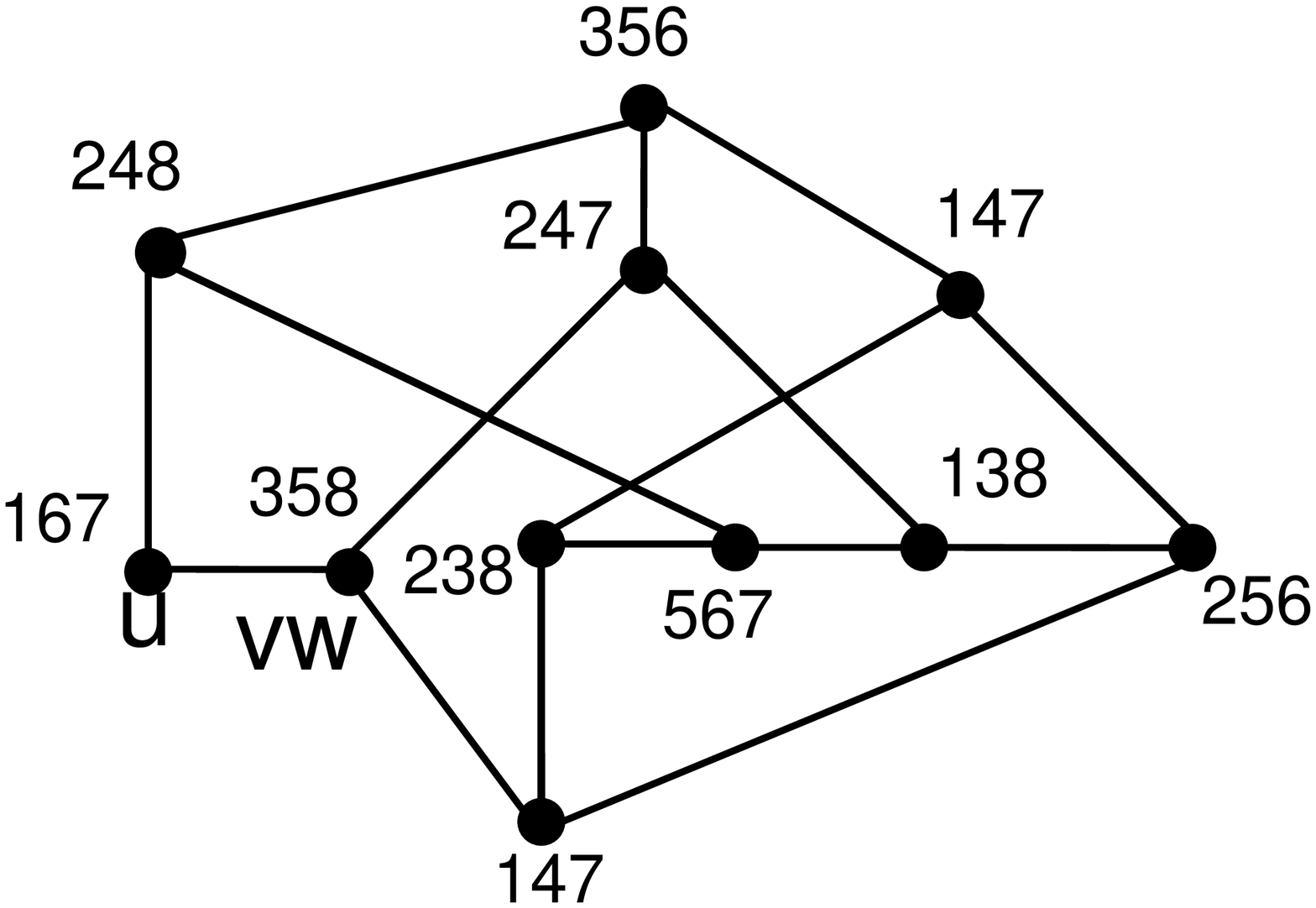, width=0.4\textwidth}}
    \hfil \psfig{figure=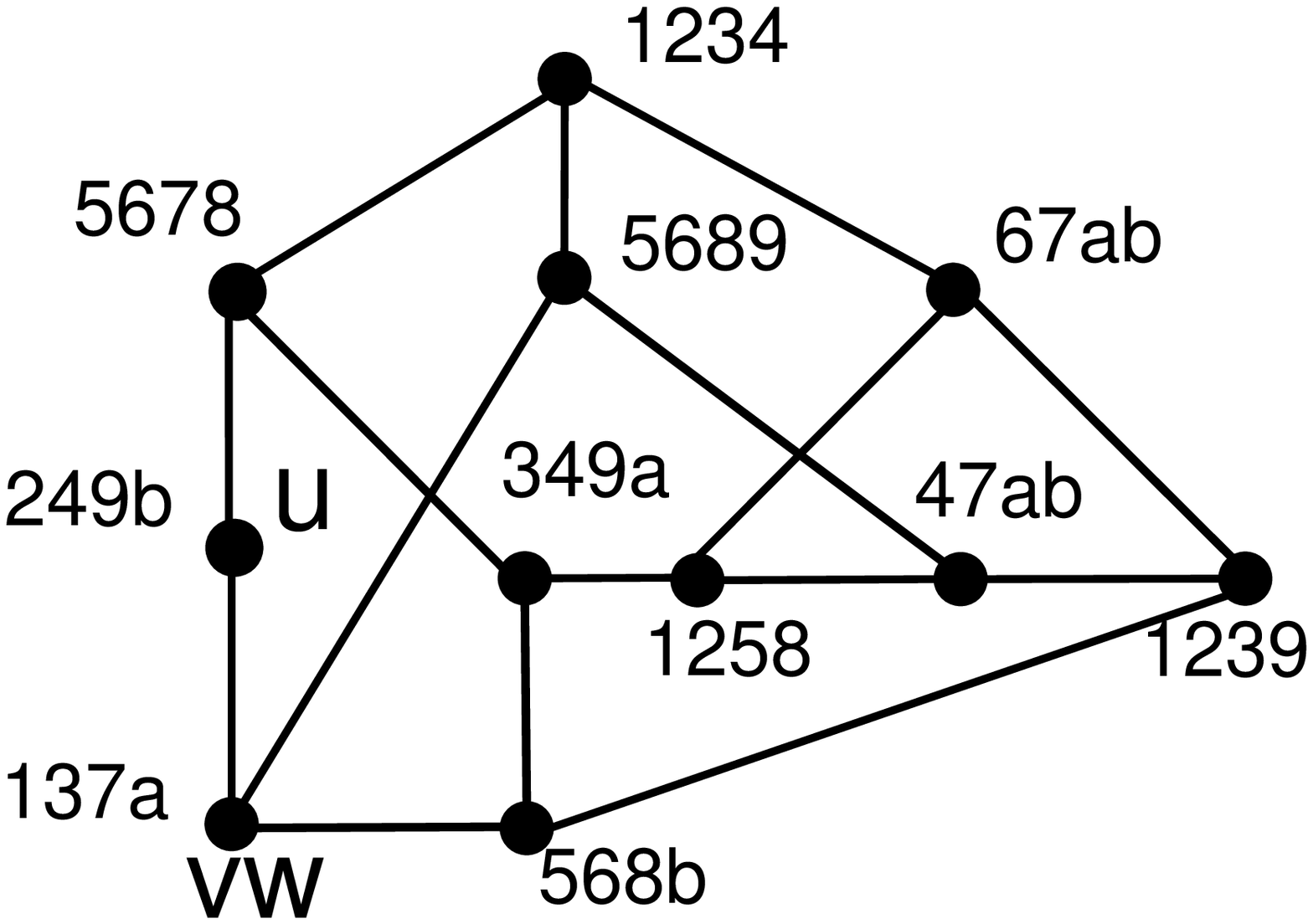, width=0.4\textwidth}}
\centerline{$(G_2/vw)+ u(vw)$ in (IV) \hspace*{1cm} \hfil \hspace*{1cm}$(G_2/vw)+ u(vw)$ in (V)}
\caption{$(G_2/vw)+ u(vw)$ in (IV) is  $8\!:\!3$-colorable,
while $(G_2/vw)+ u(vw)$ in (V) is  $11\!:\!4$-colorable.}
\label{fig:u-vw}
\end{figure}

\item[(c)] Let $(G_2/uw)+ v(uw)$ be the graph obtained by
identifying $u$  and $w$ as one vertex $uw$ followed by adding an
edge $v(uw)$. The fractional coloring $(0,1,0,0)$ is extensible in
$\F_{t_0}(G_2)$ if and only if $\chi_f((G_2/uw)+ v(uw))\leq t_0$,
which is verified by Figure \ref{fig:v-uw}.
\begin{figure}[htbp]
  \centerline{ {\psfig{figure=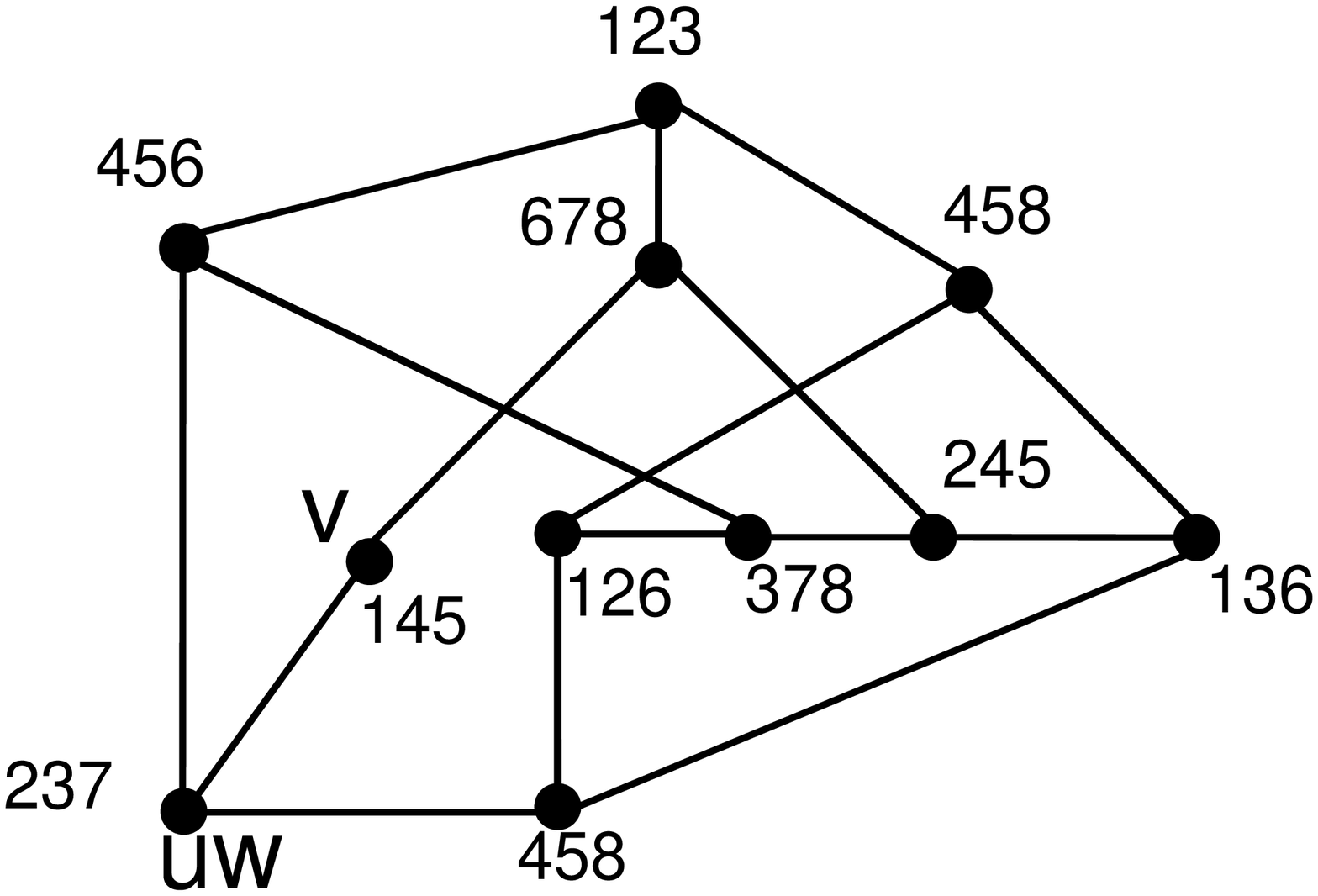, width=0.4\textwidth}}
    \hfil \psfig{figure=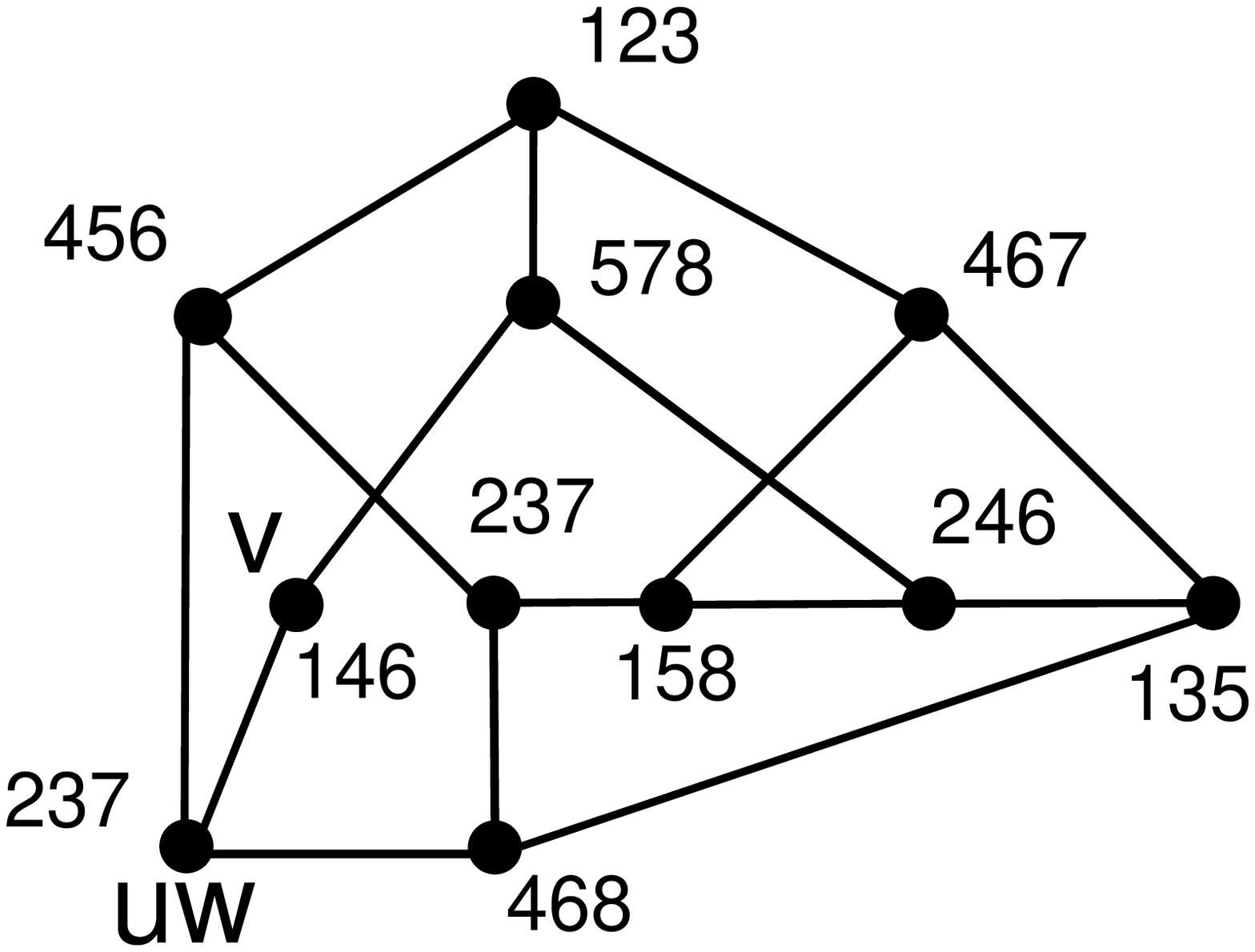, width=0.4\textwidth}}
\centerline{$(G_2/uw)+ v(uw)$ in (IV) \hspace*{1cm} \hfil \hspace*{1cm}$(G_2/uw)+ v(uw)$ in (IV)}
\caption{Both $(G_2/uw)+ v(uw)$ in (IV) and (V) are  $8\!:\!3$-colorable. }
\label{fig:v-uw}
\end{figure}

\item[(d)] Let $(G_2/uv)+ w(uv)$ be the graph obtained by
identifying $u$ and $v$ as one vertex $uv$ followed by  adding an
edge $w(uv)$. The fractional coloring $(1,0,0,0)$ is extensible in
$\F_{t_0}(G_2)$ if and only if $\chi_f((G_2/uv)+ w(uv))\leq t_0$,
which is verified by Figure \ref{fig:w-uv}.
\begin{figure}[htbp]
  \centerline{ {\psfig{figure=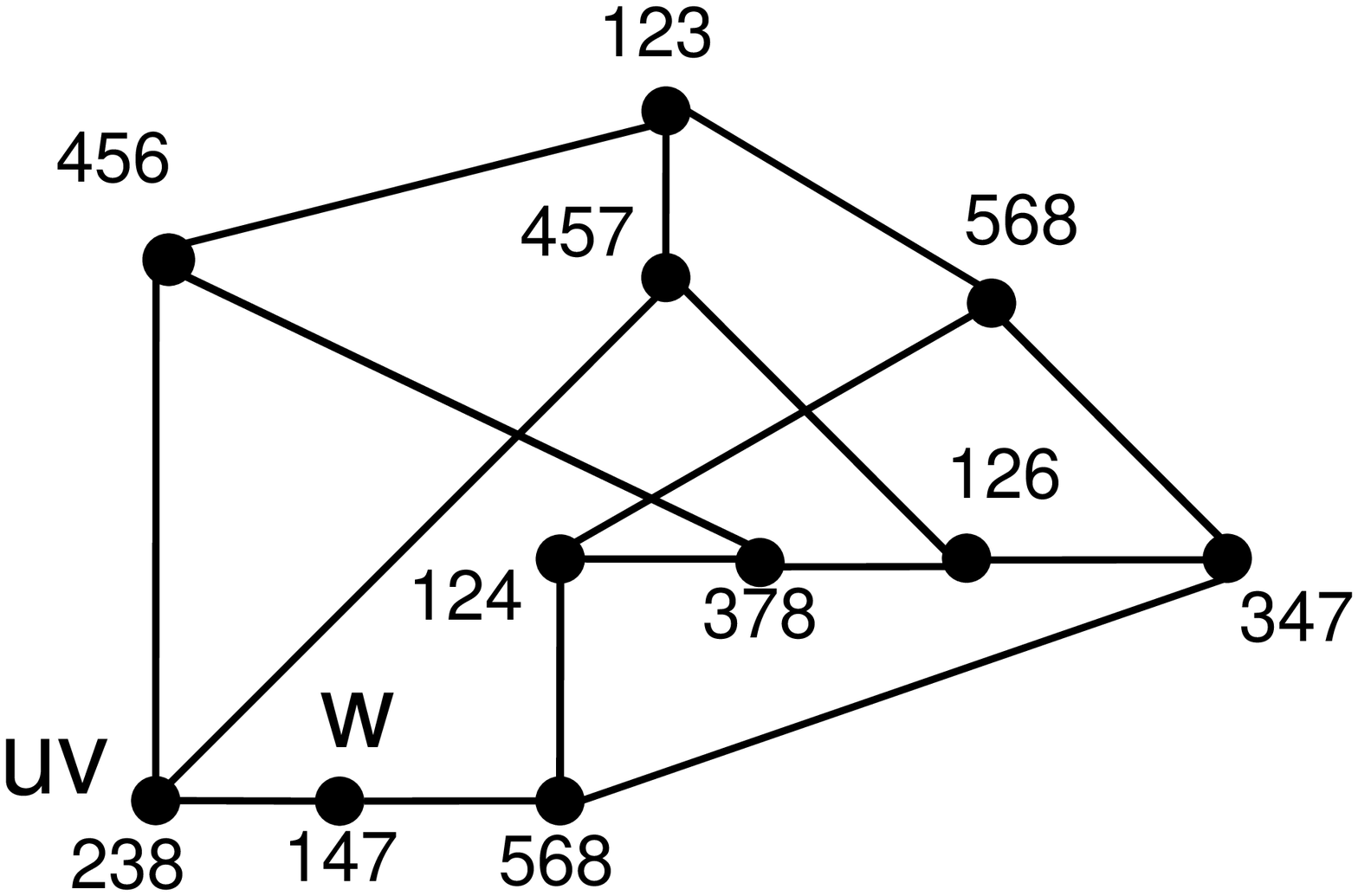, width=0.4\textwidth}}
    \hfil \psfig{figure=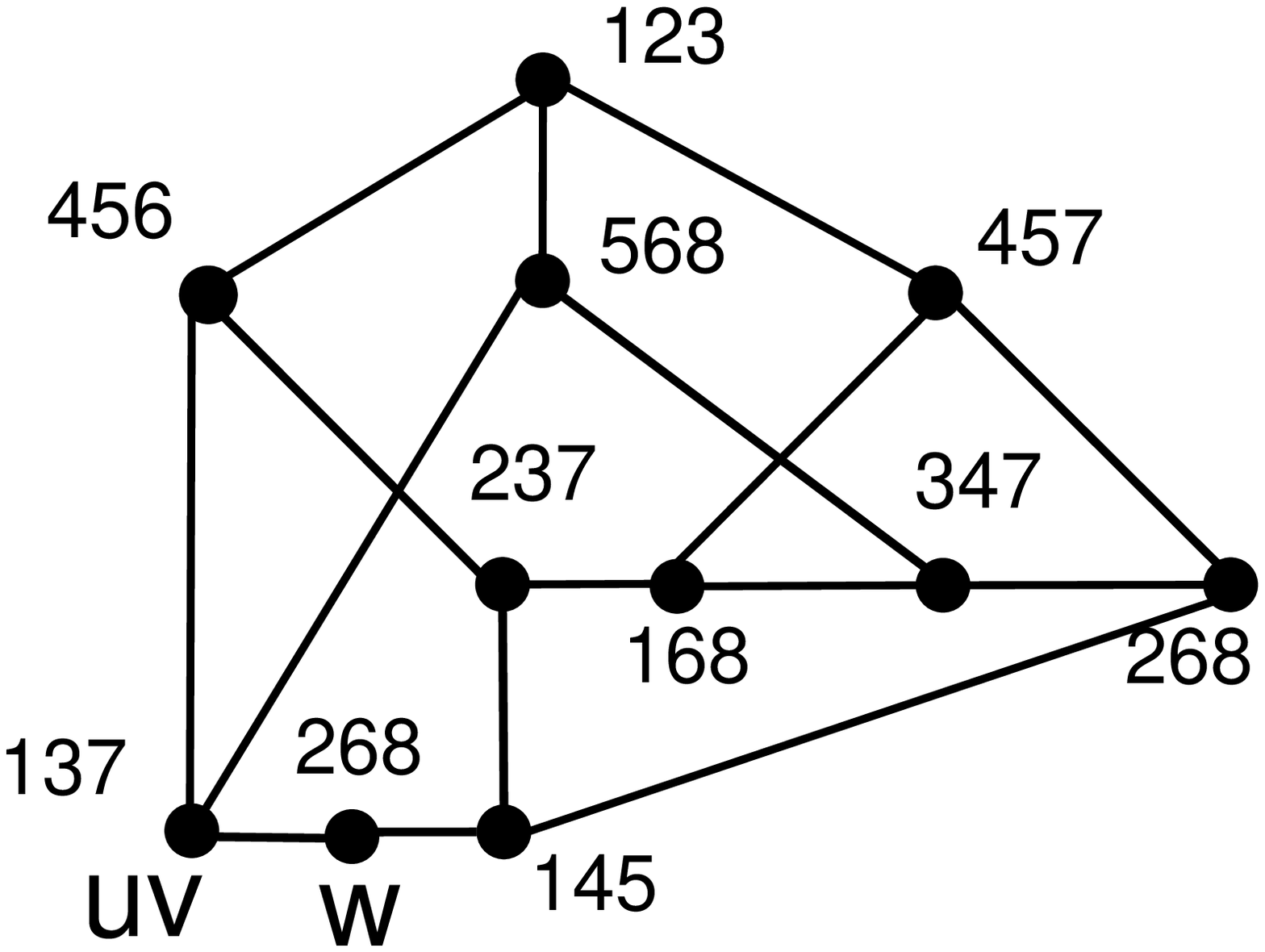, width=0.4\textwidth}}
\centerline{$(G_2/uv)+ w(uv)$ in (IV) \hspace*{1cm} \hfil
\hspace*{1cm} $(G_2/uv)+ w(uv)$ in (V)} \caption{Both $(G_2/uv)+
w(uv)$ in (IV) and (V) are  $8\!:\!3$-colorable. } \label{fig:w-uv}
\end{figure}

\item[(e)] Choose $\lambda=3t_0-8$. Since $\frac{8}{3}< \frac{11}{4}<t_0<3$, we have
$0<\lambda<1$. Note that
$$(3-t_0,0,0,0)=\lambda (0,0,0,0) + (1-\lambda)(\tfrac{1}{3},0,0,0).$$
Observe that $(3-t_0,0,0,0)$ being fully extensible in $\F(G_2)$
implies that $(3-t_0,0,0,0)$ is extensible in $\F_{t_0}(G_2)$.  To
show $(3-t_0,0,0,0)$ is fully extensible in $\F(G_2)$, it suffices
to show both $(0,0,0,0)$ and $(\frac{1}{3},0,0,0)$ are fully
extensible in $\F(G_2)$ by Lemma \ref{fully} (see Figure
\ref{fig:f0} and \ref{fig:f1}).
\begin{figure}[htbp]
  \centerline{ {\psfig{figure=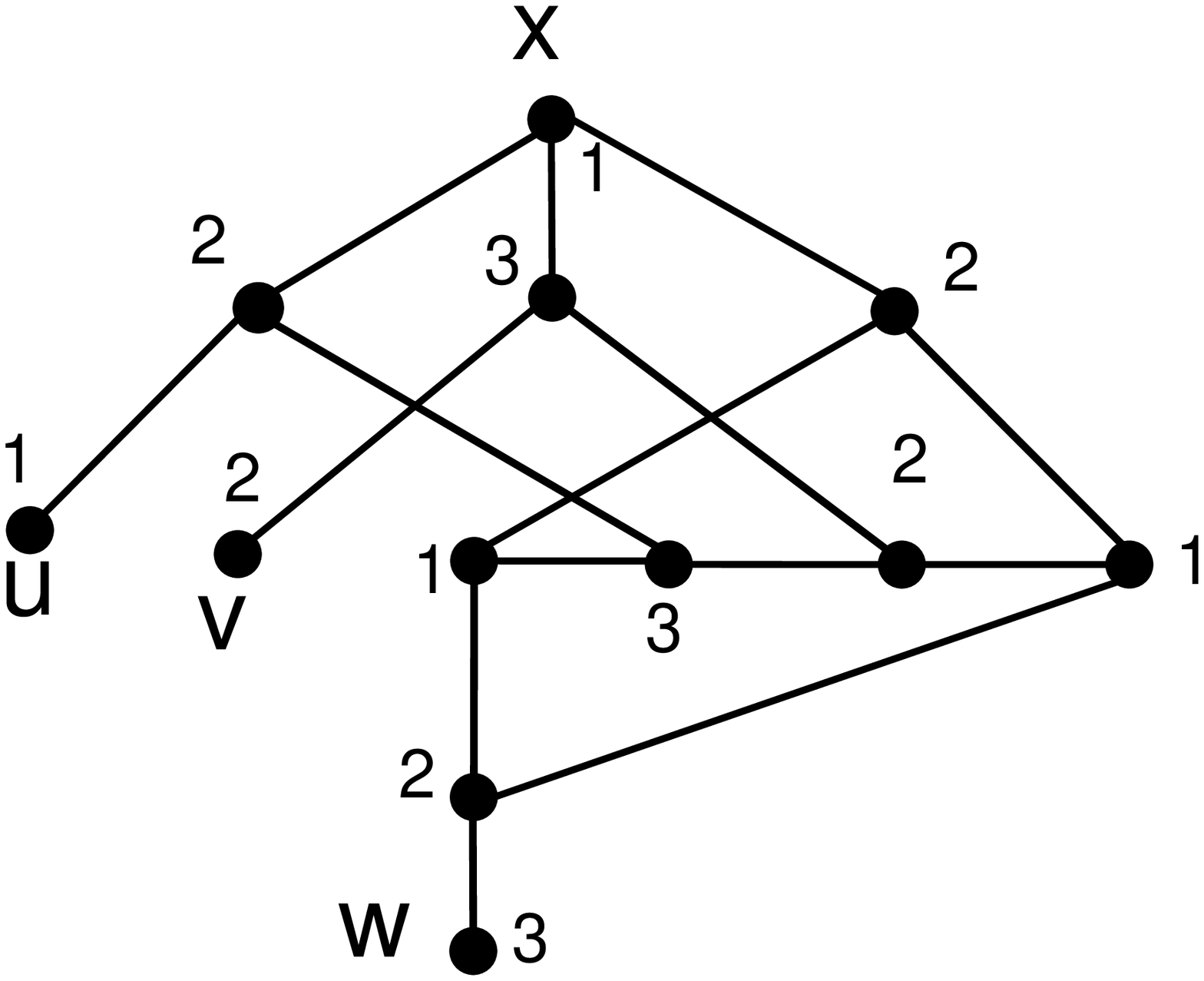, width=0.4\textwidth}}
    \hfil \psfig{figure=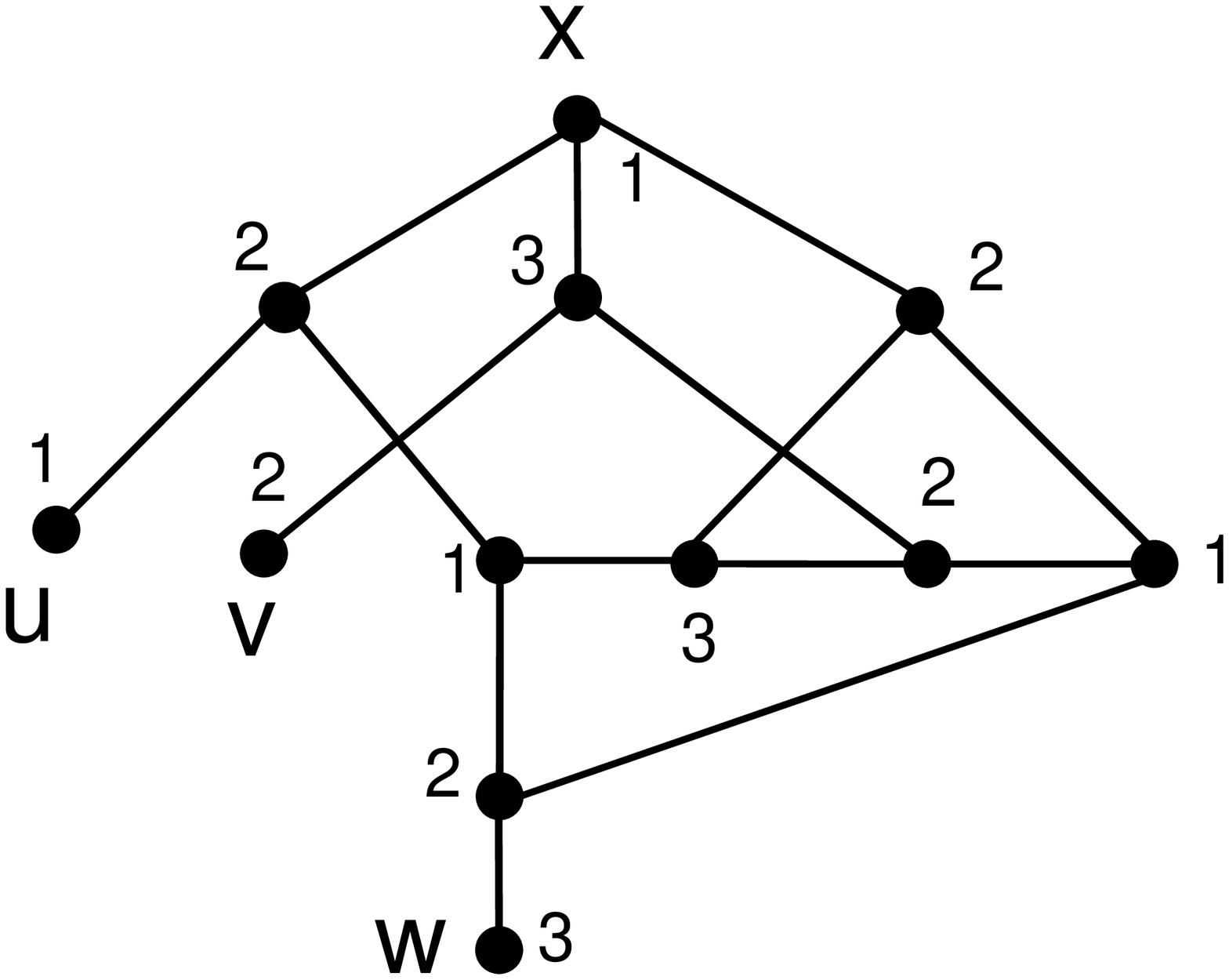, width=0.4\textwidth}}
\centerline{ (IV) \hspace*{1cm}\hfil \hspace*{1cm}(V)}
\caption{The fractional coloring $(0,0,0,0)$ is fully extensible in
$\F(G_2)$.}
\label{fig:f0}
\end{figure}

\begin{figure}[htbp]
  \centerline{ {\psfig{figure=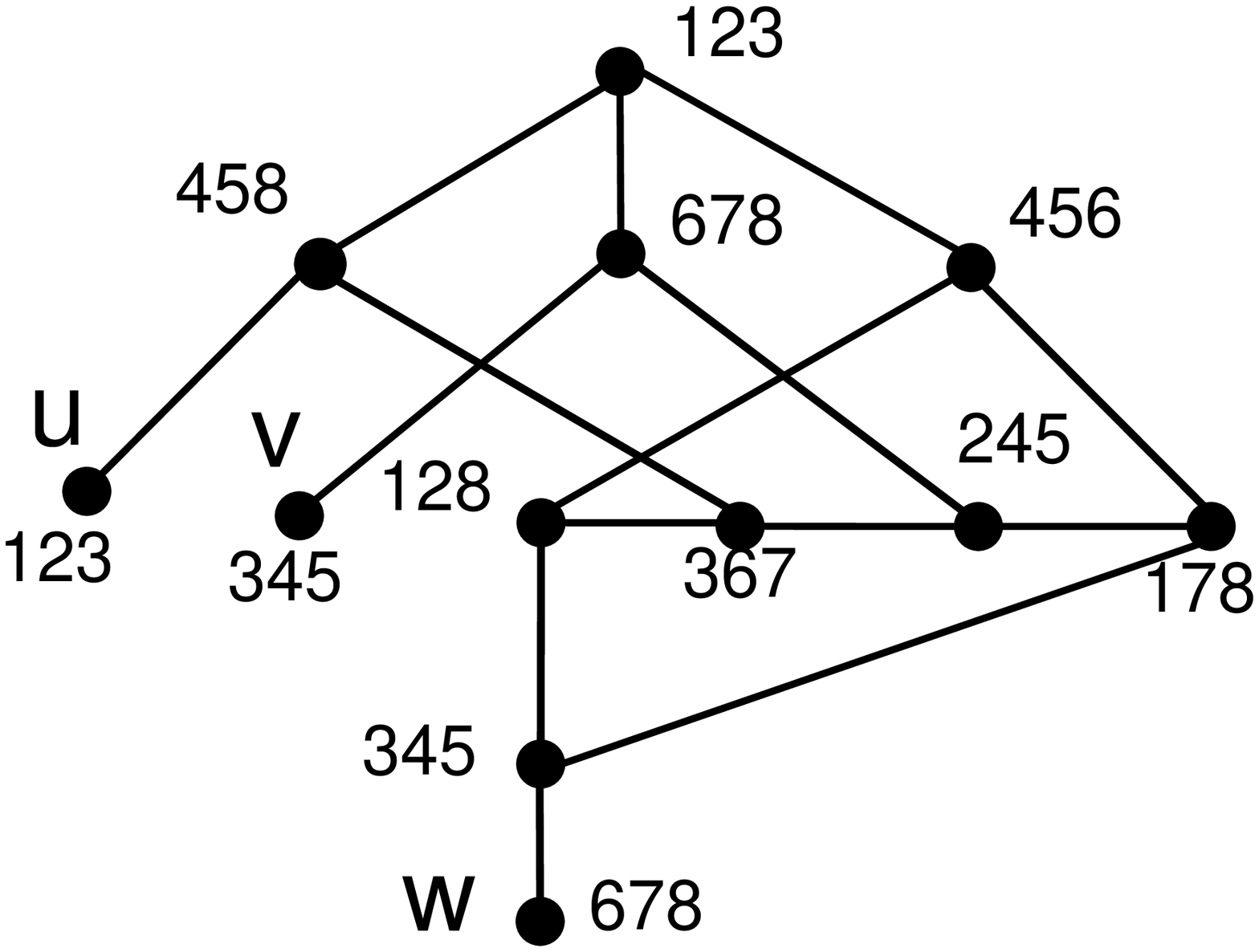, width=0.4\textwidth}}
    \hfil \psfig{figure=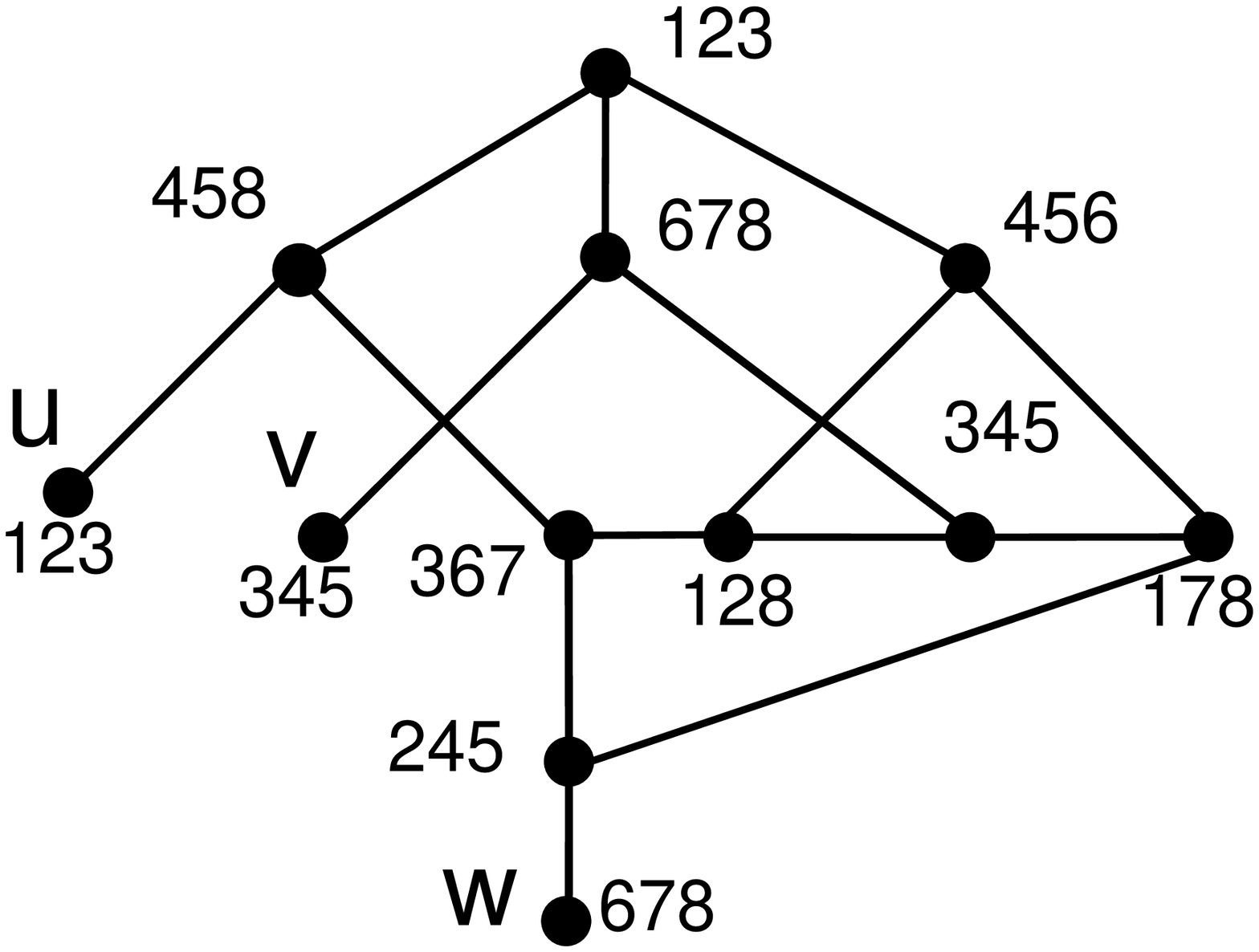, width=0.4\textwidth}}
\centerline{ (IV) \hspace*{1cm}\hfil \hspace*{1cm}(V)}
\caption{The fractional coloring $(\frac{1}{3},0,0,0)$ is fully extensible in
$\F(G_2)$.}
\label{fig:f1}
\end{figure}

\item[(f)]Similarly, to show $(0,3-t_0,0,0)$ is extensible in
  $\F_{t_0}(G_2)$, it suffices to show $(0,0,0,0)$ and $(0,
  \frac{1}{3},0,0)$ are fully extensible in $\F(G_2)$ (see Figure
  \ref{fig:f0} and \ref{fig:f2}).
\begin{figure}[htbp]
  \centerline{ {\psfig{figure=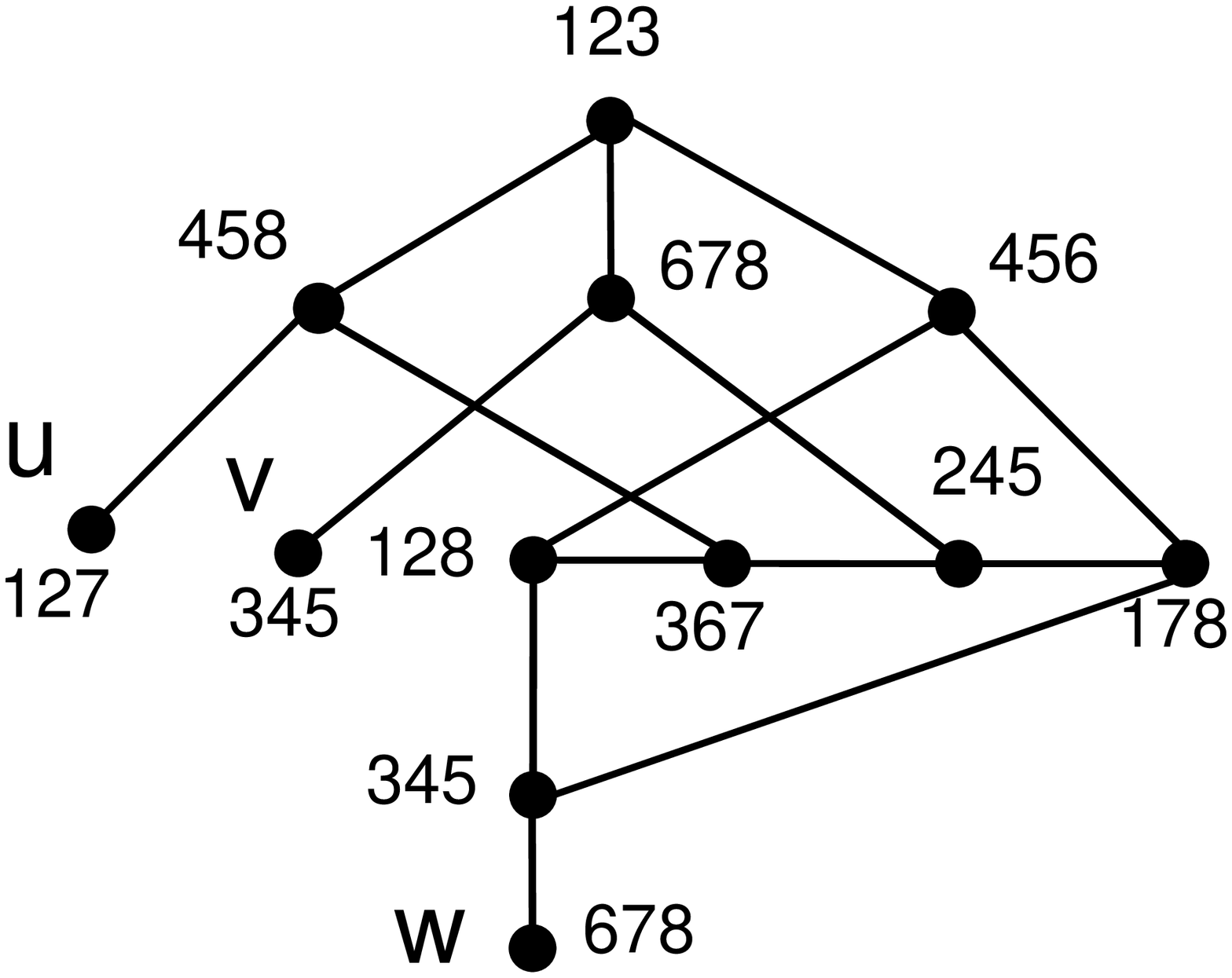, width=0.4\textwidth}}
    \hfil \psfig{figure=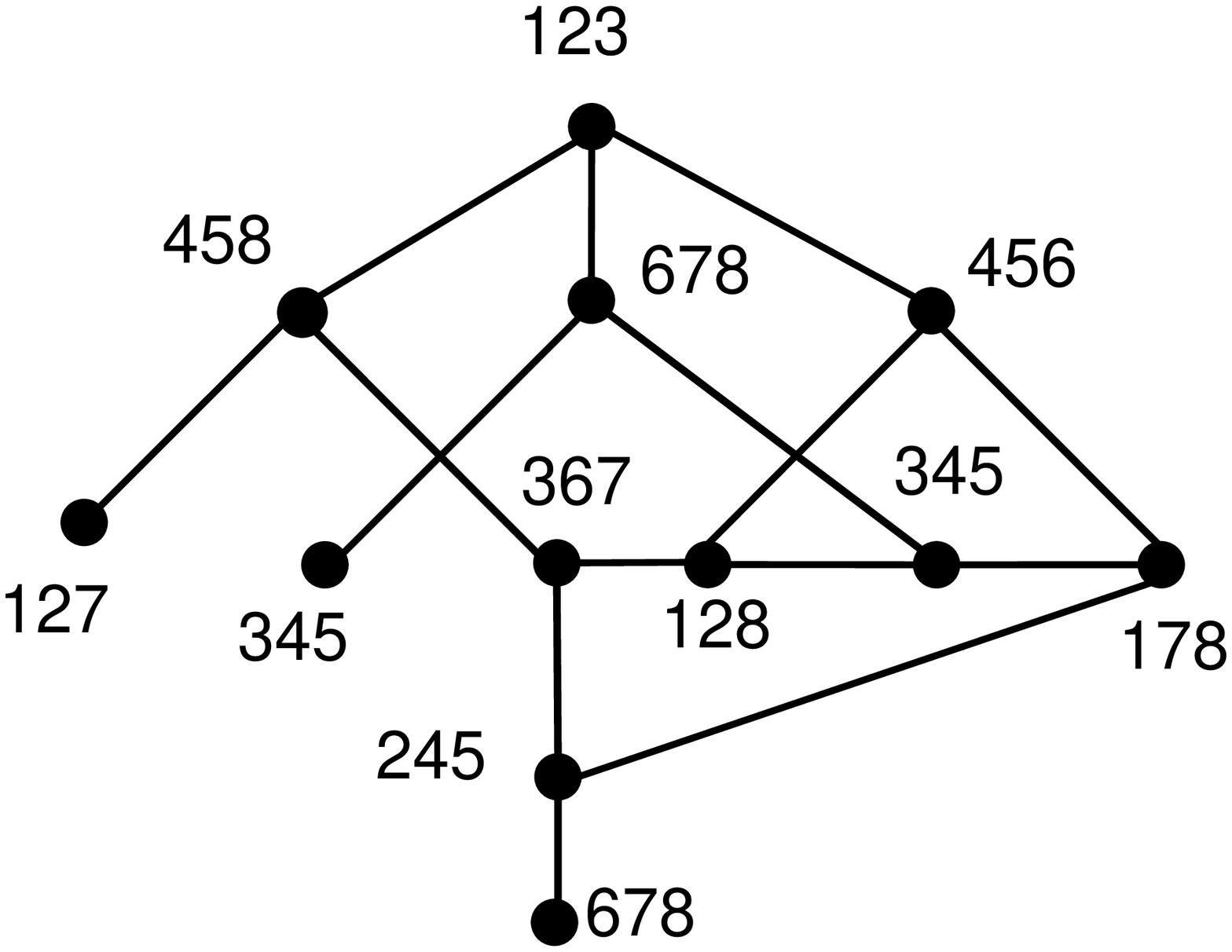, width=0.4\textwidth}}
\centerline{ (IV) \hspace*{1cm}\hfil \hspace*{1cm}(V)}
\caption{The fractional coloring $(0,\frac{1}{3},0,0)$ is fully extensible in
$\F(G_2)$.}
\label{fig:f2}
\end{figure}
\item[(g)]Similarly, to show $(0, 0,3-t_0,0)$ is extensible in
  $\F_{t_0}(G_2)$, it suffices to show $(0,0,0,0)$ and $(0, 0,
  \frac{1}{3},0)$ are fully extensible in $\F(G_2)$ (see Figure
  \ref{fig:f0} and \ref{fig:f3}).
\begin{figure}[htbp]
  \centerline{ {\psfig{figure=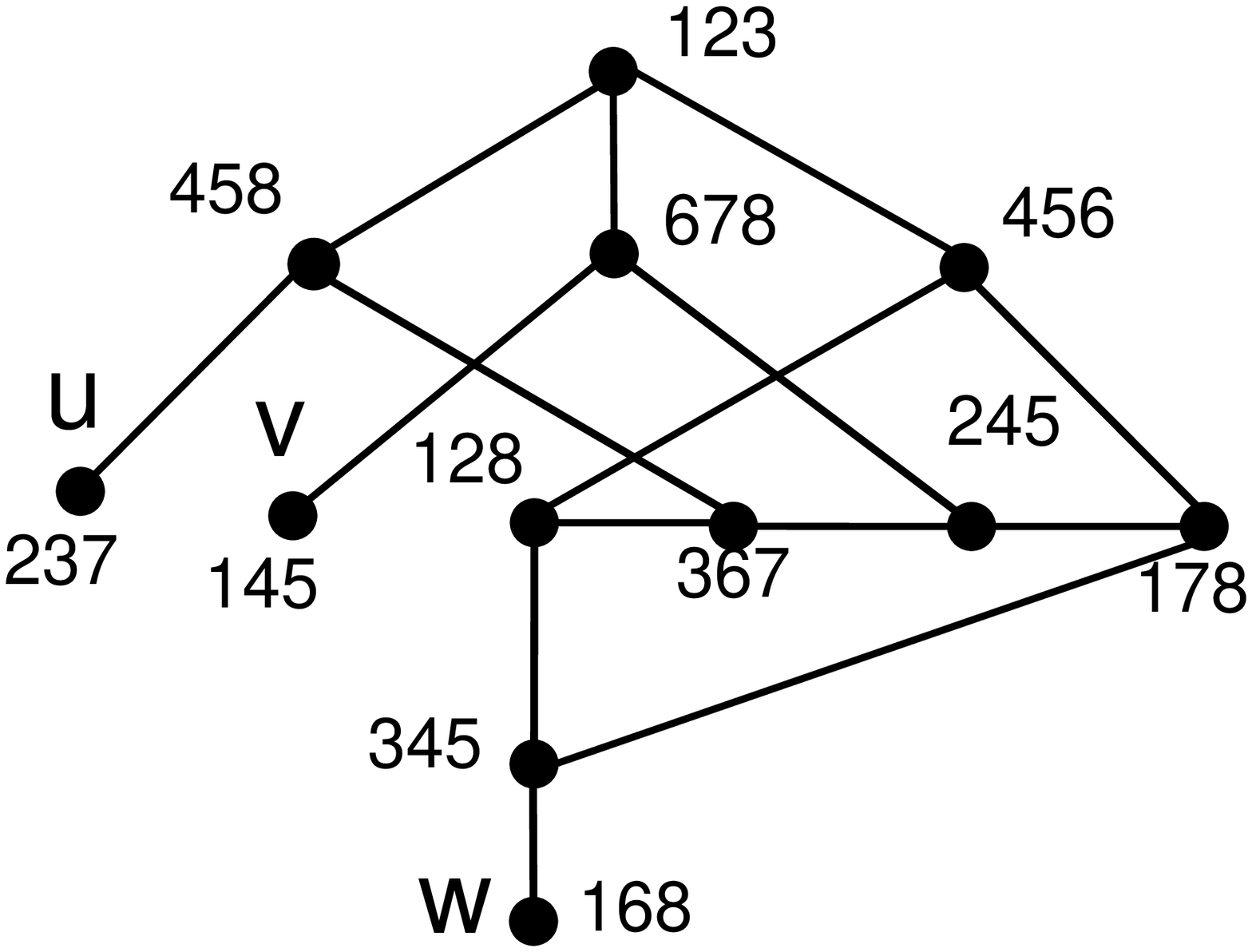, width=0.4\textwidth}}
    \hfil \psfig{figure=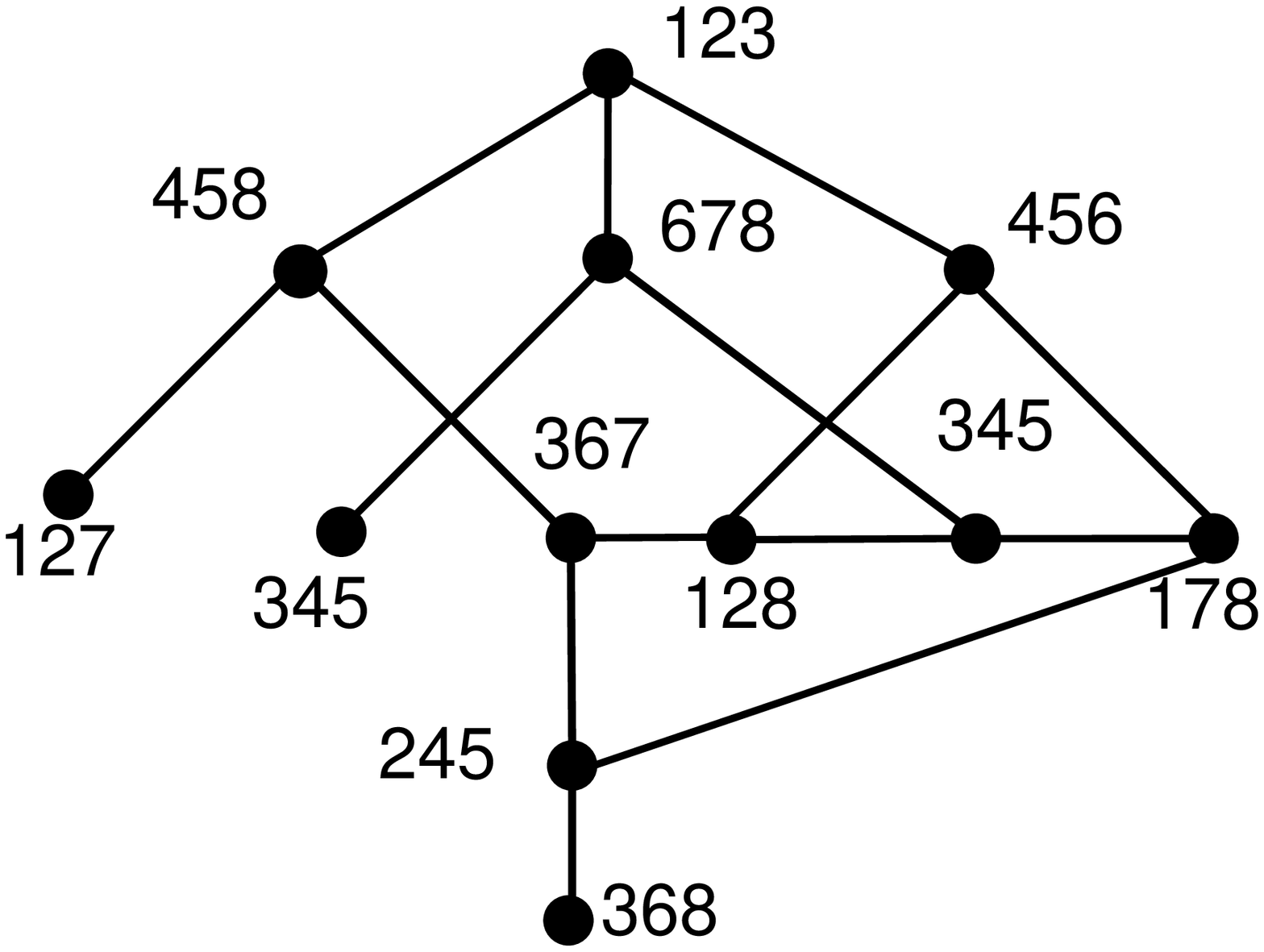, width=0.4\textwidth}}
\centerline{ (IV) \hspace*{1cm}\hfil \hspace*{1cm}(V)}
\caption{The fractional coloring $(0,0,\frac{1}{3},0)$ is fully extensible in
$\F(G_2)$.}
\label{fig:f3}
\end{figure}
\end{description}
Applying Theorem \ref{extension}, we have $\chi_f(G) \leq t_0 < t=
\chi_f(G)$,  which is a contradiction and so
$\chi_f(G_1)=t=\chi_f(G)$. However, $G$ is fractionally-critical.
Contradiction!

Now we consider the degenerated cases. For graphs (I), (II), and
(III), $\{u,v\}$ is degenerated into a set of size one. Since $G$
is $2$-connected, $G$ is a subgraph of one of graphs listed in Figure
\ref{fig:c5s2}. Thus $G$ is 8:3-colorable. Contradiction! For graph
(IV) and (IV), $\{u,v,w\}$ is degenerated into a set $H$ of size at
most $2$. If $|H|=1$, then $G$ is a subgraph of one of graphs in Figure
\ref{fig:uvw}. If $H=\{u',v'\}$, then $G_2/u'v'$ and $G_2+u'v'$ are
subgraphs of graphs from Figure \ref{fig:uvw} to Figure
\ref{fig:w-uv}. Applying Lemma \ref{l:cut2}, we get
$$\chi_f(G) \leq \max\left\{\chi(G_1), \frac{11}{4}\right \} < \chi_f(G).$$ Contradiction!
Hence, Lemma \ref{c5} follows.
\hfill $\square$

\begin{figure}[htbp]
  \centerline{ {\psfig{figure=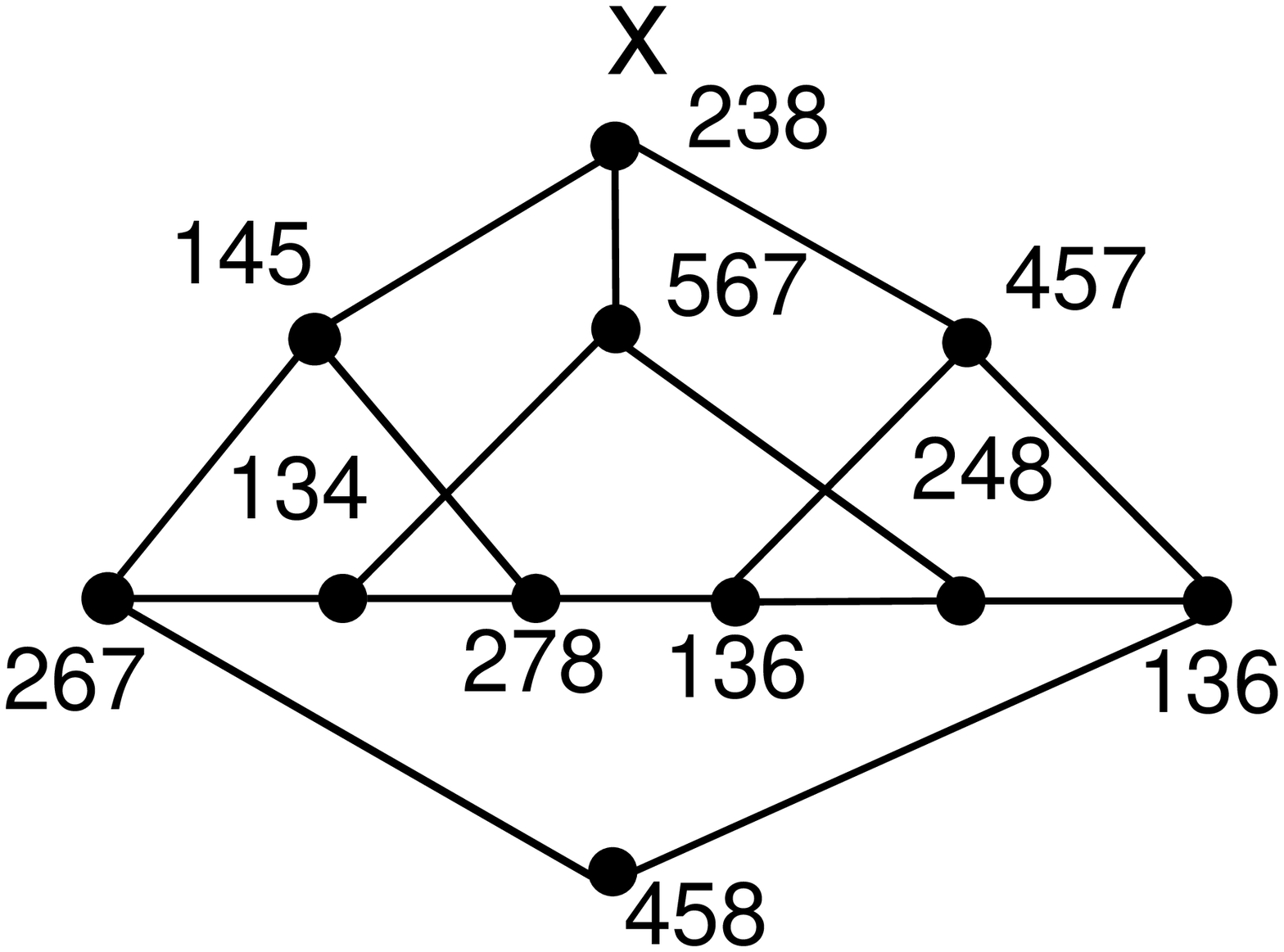, width=0.4\textwidth}}
    \hfil \psfig{figure=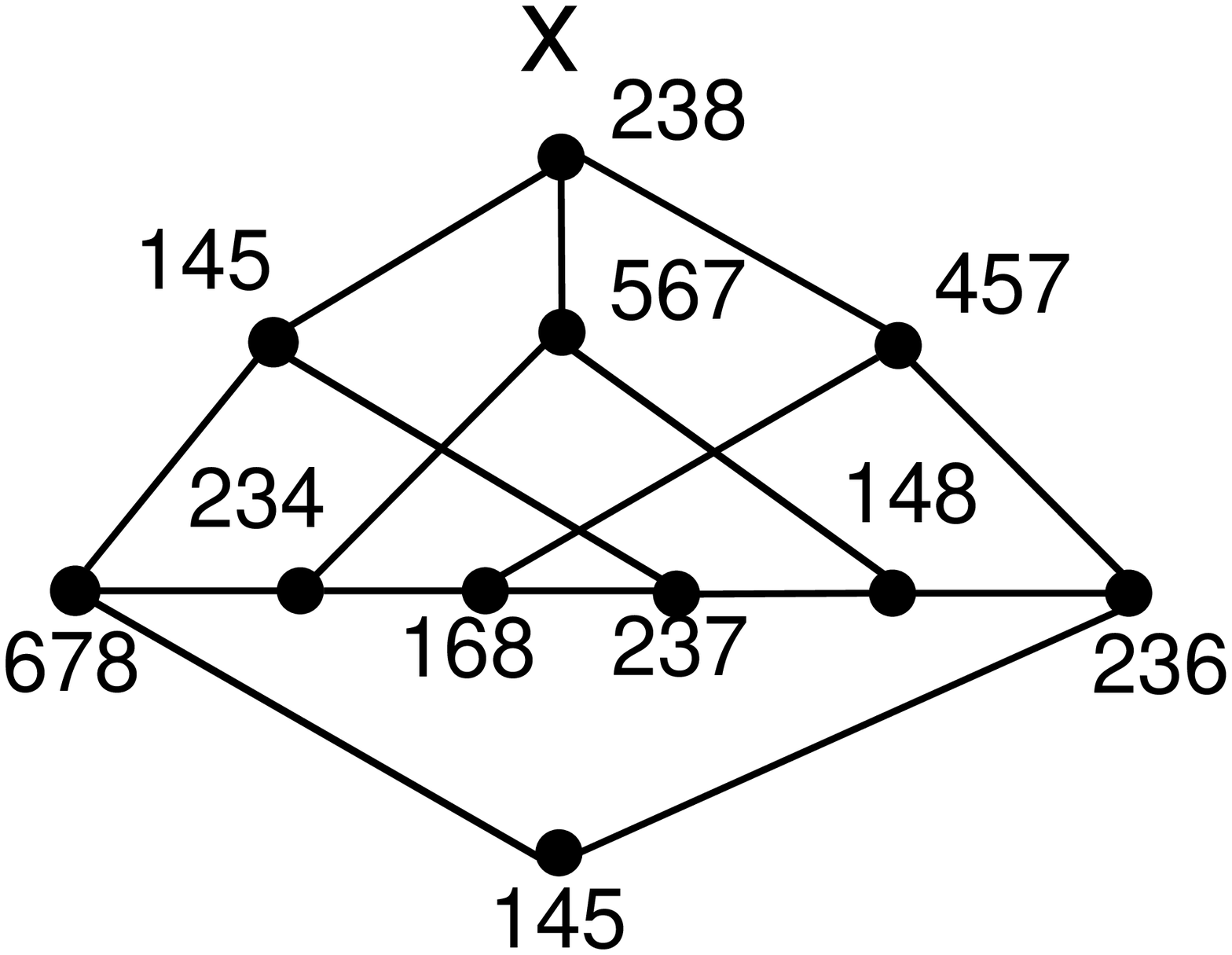, width=0.4\textwidth}}
  \centerline{ {\psfig{figure=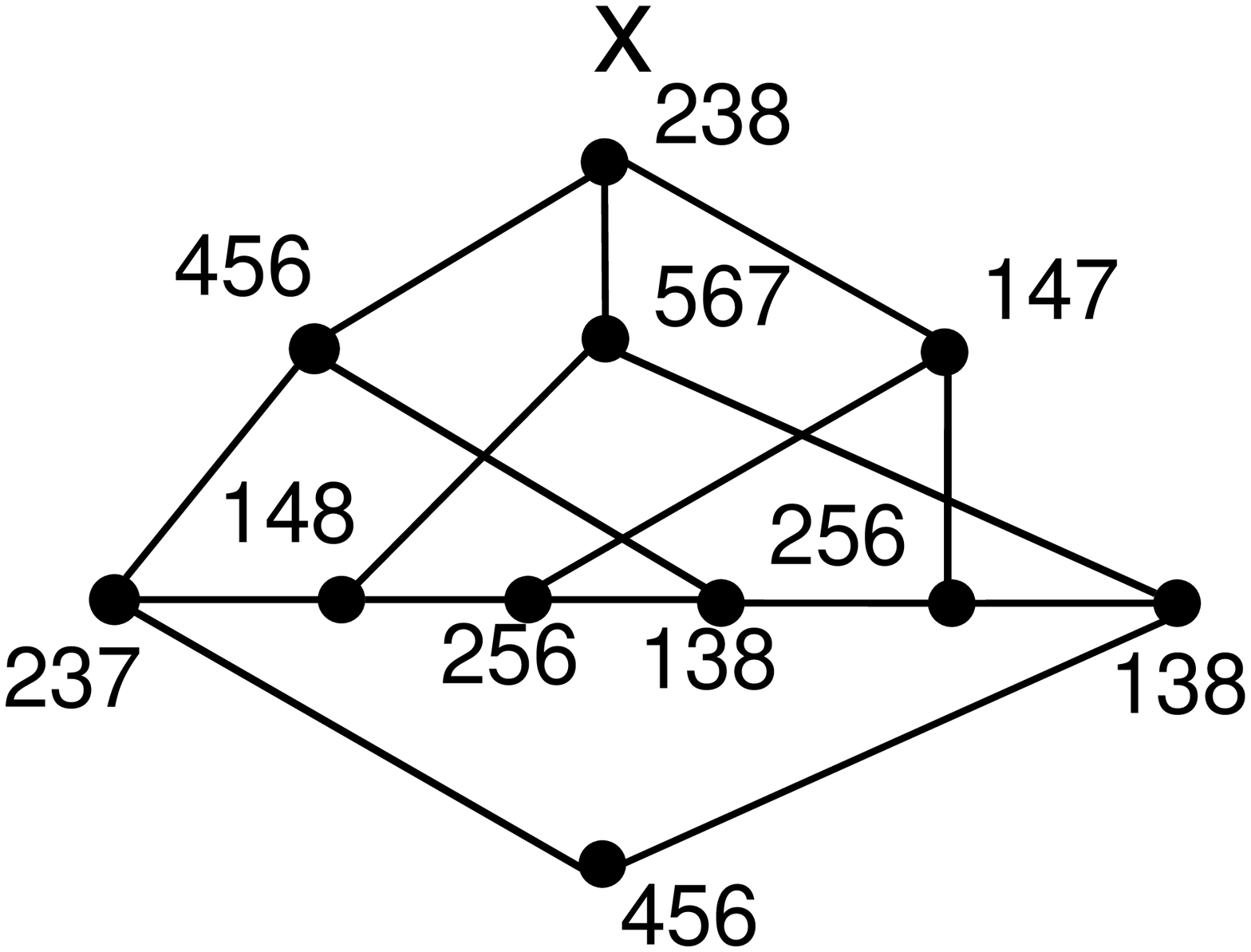, width=0.4\textwidth}}
    \hfil \psfig{figure=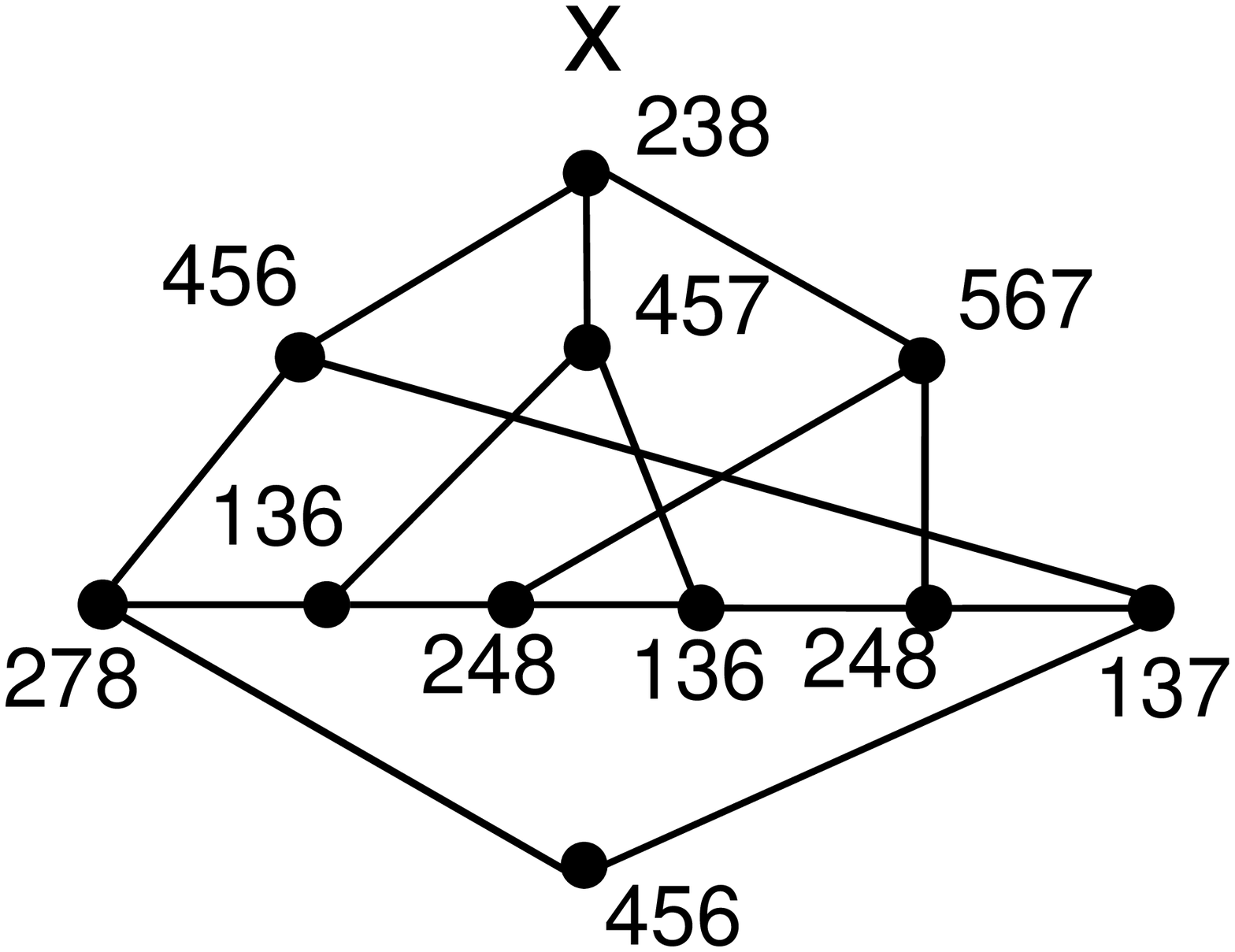, width=0.4\textwidth}}
\caption{All possible graphs with $|V(C)\cap N_G^2(x)|=6$.}
\label{fig:c7}
\end{figure}

\begin{lemma}\label{c7}
Assume that $G$ is a fractionally-critical triangle-free graph satisfying
$\Delta(G) \leq 3$ and $ \frac{8}{3} < \chi_f(G) < 3$. For any
vertex $x \in V(G)$ and any $7$-cycle $C$ of $G$, we have $|V(C)\cap N_G^2(x)|\leq
5$.
\end{lemma}
\noindent {\bf Proof:} We prove the statement by contradiction. Suppose that there
is a vertex $x$ and a $7$-cycle $C$ satisfying $|V(C)\cap
N_G^2(x)| \geq 6$. Recall that $|N_G^2(x)| \leq 6$. We have  $|V(C)\cap
N_G^2(x)| = 6$.  Combined with the
fact that $G$ is triangle-free and 2-connected, $G$ must be one of the
following graphs in Figure  \ref{fig:c7}.
All of them are  $8\!:\!3$-colorable, see Figure \ref{fig:c7}.  Contradiction!

\hfill$\square$

\section{Admissible sets and Theorem \ref{t2}}
The following approach is similar to the one used in \cite{hz}. A
significant difference is a new concept ``admissible set''.
Basically, it replaces the independent set $X$ (of $G^\ast$ in
\cite{hz}) by three independent sets $X_1\cup X_2\cup X_3$.

 Recall that $v\in N^i_G(u)$ if there exists a $uv$-path of length $i$
 in $G$. A set $X \subset V(G)$ is called  {\it admissible} if
 $X$ can be partitioned into three sets $X_1$, $X_2$, $X_3$ 
 satisfying
\begin{enumerate}
\item  If $\{u, v\} \subset X_i$ for some $i \in \{1,2,3\}$, then
$v\not\in N_G^1(u) \cup N_G^3(u) \cup N_G^5(u)$.
\item  If $u\in X_i$ and $v\in X_j$ for some $i$ and $j$ satisfying
 $1\leq i\not =j \leq 3$, then
$v\not\in  N_G^1(u) \cup N_G^2(u) \cup N_G^4(u)$.
\end{enumerate}

The following key theorem connects the $\chi_f(G)$ to a partition of
$G$ into admissible sets. 
\begin{theorem} \label{t2} Assume that  $G$ is a fractionally-critical
  triangle-free graph satisfying  $\Delta(G) \leq 3$  and
  $\frac{11}{4}<\chi_f(G) <3$.  If $V(G)$ can be partitioned
  into $k$ admissible sets, then
\begin{equation}
  \chi_f(G)\leq 3 -\frac{3}{k+1}.
\end{equation}
\end{theorem}
Let $X=X_1\cup X_2 \cup X_3$ be an admissible set.
Inspired by the method used in \cite{hz}, we define an auxiliary graph
$G'=G'(X)$ as follows. We use the notation
$\Gamma(X)$ to denote the neighborhood of $X$ in $G$.
For each $i \in \{1,2,3\}$, let $Y_i=\Gamma(X_i)$.
By admissible conditions, $Y_1$, $Y_2$, and $Y_3$ are all independent
sets of $G$.   Let $G'$ be a graph obtained from $G$ by
deleting $X$; identifying each $Y_i$ as a single vertex $y_i$ for
$1\leq i\leq 3$; and adding three edges $y_1 y_2$, $y_2 y_3$, $y_1 y_3$.
 We have the following lemma, which will be proved later.
 \begin{lemma} \label{l:3colorable} Assume that $G$ is
   fractionally-critical triangle-free graph satisfying $\Delta(G) \leq 3$ and  $\frac{11}{4} < \chi_f(G) < 3$.  Let $X$ be an admissible set
   of $G$ and $G'(X)$ be the graph defined as above. We have  $G'(X)$ is
   $3$-colorable.
\end{lemma}

\noindent {\bf Proof of Theorem \ref{t2}:} Assume $G$ can be
partitioned into $k$ admissible sets, say $V(G)=\cup_{i=1}^k X_i$,
 where $X_i=X_i^1\cup X_i^2 \cup X_i^3$. For each $1 \leq i \leq k $ and $1 \leq j
\leq 3$, let $Y_i^j=\Gamma(X_i^j)$. From the definition of an
admissible set and $G$ being triangle-free, we have $Y_i^j$ is an
independent set for all $1 \leq i \leq k $ and $1 \leq j \leq 3$.

By Lemma \ref{l:3colorable}, $G'(X_i)$ is $3$-colorable for all
$1\leq i \leq k$. Let $c_i$ be a 3-coloring of $G'(X_i)$ with the
color set $\{s_i^1, s_i^2, s_i^3\}$. Here all the colors $s_i^j$'s are pairwise
distinct. We use $\mathcal P(S)$ to
denote the set of all subsets of $S$. We define $f_i: V(G)
\rightarrow \mathcal P(\{s_i^1,s_i^2,s_i^3\})$ satisfying
 $$
 f_i(v)= \left\{
 \begin{array}{ll}
 \{c_i(v)\}  \   & \textrm{if} \  v \in V(G)-  (\cup_{j=1}^3 X_i^j \cup _{j=1}^3 Y_i^j), \\
 \{c_i(y_
 i^j)\}   \  &  \textrm{if} \  v \in Y_i^j,\\
 \{s_i^1,s_i^2,s_i^3\} - c_i(y_i^j) \   &  \textrm{if} \  v \in
 X_i^j
.
 \end{array}
 \right.
 $$
 Note that for a fixed $1 \leq i \leq k$, $y_i^j$ denotes the vertex
 of $G'(X_i)$ obtained from contracting $Y_i^j$ for $1 \leq j \leq 3$.
 Observe that each vertex in $X_i^j$ receives two colors from $f_i$
 and every other vertex receives one color.  It is clear that any two
 adjacent vertices receive disjoint colors. Let $\sigma: V(G)
 \rightarrow \mathcal P(\cup_{i=1}^{k} \{s_i^1,s_i^2,s_i^3\})$ be a
 mapping defined as $\sigma(v)=\cup_{i=1}^{k} f_i(v)$. Now $\sigma$ is
 a $(k+1)$-fold coloring of $G$ such that each color is drawn from a
 palette of $3k$ colors. Thus we have  $\chi_f(G) \leq
 \frac{3k}{k+1}=3-\frac{3}{k+1}$.

 We completed the proof of  theorem \ref{t2}. \hfill $\square$

Before we prove Lemma \ref{l:3colorable}, we first prove a lemma on coloring the
 graph obtained by splitting the hub of an odd wheel. Let
 $\{x_0,x_1,\ldots,x_{2k}\}$ be the set of vertices of an odd cycle $C_{2k+1}$
 in a circular order. Let $Y=\{y_1,y_2,y_3\}$. We construct a graph $H$ as
 follows:
\begin{enumerate}
\item  $V(H)=V(C) \cup Y $.

\item $E(C) \subset E(H)$.

\item  Each $x_i$ is adjacent to exactly one element of $Y$.

\item $y_1y_2, y_2y_3, y_1y_3\in E(H)$.

\item $H$ can have at most one vertex (of $y_1$, $y_2$, and $y_3$) with degree
$2$.
\end{enumerate}

The graph $H$ can be viewed as splitting the hub of the odd wheel
into $3$ new hubs where each spoke has to choose one new hub to
connect, then connect all new hubs. One special case it that one of
the new hubs has no neighbor in $C$.

\begin{lemma}
\label{lemma1} A graph $H$ constructed as described above is $3$-colorable.
\end{lemma}
\noindent {\bf Proof:} Without loss of generality, we assume $d_H(y_1)\geq 3$
and $d_H(y_2)\geq 3$. We construct a proper $3$-coloring $c$ of $H$
as follows. First, let $c(y_1)=1$, $c(y_2)=2$, and $c(y_3)=3$.

Again, without loss of generality, we assume $(x_0,y_1) \in E(H)$.
The neighbors of $y_1$ divide $V(C)$ into several intervals. The
vertices in each interval are either connecting to $y_2$ or $y_3$
but not connected to $y_1$. As $v$ goes through each interval counter-clockwisely,
we list the neighbor of $v$ (in $Y$) and get a sequence consisting of $y_2$ $y_3$.
Then we delete the repetitions of $y_2y_2$ and
$y_3y_3$ in the sequence. There are $4$ types of intervals
based on the result seqeunce:

Type I: $y_2, y_3, y_2, y_3, \ldots,  y_3, y_2$.

Type II: $y_2, y_3, y_2, y_3, \ldots, y_2,  y_3$.

Type III: $y_3, y_2, y_3, y_2, \ldots,  y_3, y_2$.

Type IV: $y_3, y_2, y_3, y_2, \ldots, y_2, y_3$.

If $y_3$ is a vertex of degree two in $H$, then the interval $I$ has only
one type, which is degenerated into $y_2$.

Given an interval $I$, let $u(I)$ (or $v(I)$) be the common neighbor
of $y_1$ and the left (or right) end of $I$ respectively. We color
$u(I)$ and $v(I)$ first and then try to extend it as a proper
coloring of $I$. Sometimes we succeed while sometimes we fail. We
ask a  question whether we can always get a proper coloring.  The
answer depends  only on the type of $I$ and the coloring combination
of $u(I)$ and $v(I)$.  In Table \ref{tab:1}, the column is
classified by the coloring combination of $u(I)$ and $v(I)$, while
the row is classified by the types of $I$.  Here ``yes'' means the
coloring process always succeeds,  while ``no'' means it sometimes
fails.

\begin{table}[hbt]
  \centering
  \begin{tabular}{|c|c|c|c|c|}
    \hline
   & $(2,2)$ & $(2,3)$ & $(3,2)$ & $(3,3)$\\
\hline
Type I & Yes & Yes & Yes & No\\
\hline
Type II & Yes & Yes & No & Yes\\
\hline
Type III & Yes & No & Yes & Yes\\
\hline
Type IV & No & Yes & Yes & Yes\\
\hline
  \end{tabular}
  \caption{Can a coloring be extended to $I$ properly?}
  \label{tab:1}
\end{table}
An ending vertex $w$ of an interval $I$ is called a {\it free} end
if $w$'s two neighbors outside $I$ receiving the same color. Observe
that in all yes entries, there exist at least one free end. Note
that  each vertex on $I$ has degree $3$. We can always color the
vertices of $I$ greedily starting from the  end not equaling $w$.
Since $w$ is a free end, there is no difficulty to color $w$ at the
end.

Now we put them together. We color the neighbors of $y_1$ one by one
counter-clockwise starting from $x_0$ according to the following
rules:

\begin{enumerate}
\item When we meet an interval $I$ of type II or III, we keep the colors
of $u(I)$ and $v(I)$ the same.
\item When we meet an interval $I$ of type I or IV, we keep the colors
of $u(I)$ and $v(I)$ different.
\end{enumerate}

There are two possibilities. If the last interval obeys the rules,
then by Table \ref{tab:1}, we can extend the partial coloring into a
proper 3-coloring of $H$. If the last interval does not obey the
rules, then we swap the colors $2$ and $3$ of the neighbors of
$y_1$. By Table \ref{tab:1}, the new partial coloring can be
extended into a proper 3-coloring of $H$.
We completed the proof. \hfill $\square$

A maximal  $2$-connected subgraph $B$ of a graph is called a {\it
block} of $G$.  A {\it Gallai tree} is a connected graph in which
all blocks are either complete graphs or odd cycles. A {\it
  Gallai forest} is a graph all of whose components are Gallai
trees. A $k$-{\it Gallai tree (forest)} is a Gallai tree (forest) such
that the degree of all vertices are at most $k-1$. A {\it $k$-critical
  graph} is a graph $G$ whose chromatic number is $k$ and the
chromatic number of any proper subgraph is strictly less than $k$. Gallai
showed the following Lemma.

\begin{lemma}[Gallai \cite{gallai}]\label{lemma3}
If $G$ is a $k$-critical graph, then the subgraph of $G$ induced on
the vertices of degree $k-1$ is a k-Gallai forest.\end{lemma}
Now, we are ready to prove Lemma \ref{l:3colorable}.

\noindent {\bf Proof of Lemma \ref{l:3colorable}:} Write $G'=G'(X)$
for short. Note that  the only possible vertices of degree greater
than $3$ in $G'$ are $y_1$, $y_2$, and $y_3$. We can color $y_1$,
$y_2$, and $y_3$ by $1$, $2$, and $3$, respectively. Since the
remaining vertices have degree at most $3$, we can color $G'$
properly with $4$ colors greedily, i.e., $\chi(G') \leq 4$.

Suppose that $G'$ is not $3$-colorable. Let $H$ be a 4-critical
subgraph of $G'$. We have  $d_{H}(v) = 3 $ for all $v \in H$ except
for possible $y_1$, $y_2$, and $y_3$. Let $T$ be the subgraph
induced by all $v \in H$ such that $d_{H}(v) = 3 $; then $T$ is not
empty since $|T|\geq |H|-3\geq 1$. By Lemma \ref{lemma3} the
subgraph of $H$ induced by $T$ is a $4$-Gallai forest. We known $T$ may contain one
or more vertices in $\{y_1, y_2,y_3\}$. Let $T'=T\setminus \{y_1,
y_2,y_3\}=V(H)\setminus \{y_1, y_2,y_3\}$. Observe that any induced
subgraph of a $4$-Gallai forest is still a $4$-Gallai forest and so
the subgraph of $H$ induced by $T'$ is also a $4$-Gallai forest.

Recall the definition of an admissible set.  If $u\in X_i$ and $v\in
X_j$ for some $i$ and $j$ satisfying $1\leq i\not =j \leq 3$, then
$v\not\in N_G^1(u) \cup N_G^2(u) \cup N_G^4(u)$, which implies that
any vertex $x$ in $T'$ can have at most one neighbor in $\{y_1,
y_2,y_3\}$. Note $d_H(x)=3$. We have $d_{T'}(x)\geq 2$.

Let $B$ be a leaf block in the Gallai-forest $T'$; then $B$ is a
complete graph or an odd cycle by the definition of the
Gallai-forest. Observe  $B$ can not be a single vertex or $K_2$
since every vertex in $B$ has at least two neighbors in $T'$. As $G$
is triangle-free, then $B$ must be an odd cycle $C_{2r+1}$ with
$r\geq 2$.

\begin{description}
\item[Case (a):] $\left|N_{H}(B) \cap \{y_1,y_2,y_3\} \right| \geq 2$.
Since $H$ is $4$-critical, $H\setminus B$ is $3$-colorable. Let $c$
be a proper $3$-coloring of $H\setminus B$.  Since $H\setminus B$
contains a triangle $y_1y_2y_3$,  we have   $y_1$, $y_2$, and $y_3$
receive different colors. We have
 \[\left|N_{H}(B) \cap \{y_1, y_2,
y_3\} \right| \geq 2.\] 
By Lemma \ref{lemma1}, we can extend the
coloring $c$ to all vertices on $B$ as well. Thus $H$ is
$3$-colorable. Contradiction.

\item [Case (b):] $\left|N_{H}(B) \cap \{y_1,y_2,y_3\} \right| =
  1$.  Since $B$ is a leaf block, there is at most one vertex, say $v_0$, which is 
  connected to another block in $T'$. List the vertices of $B$ in a
  circular order as $v_0, v_1, v_2,\ldots, v_{2r}$. All
  $v_1,v_2,\ldots, v_{2r}$ connect to one $y_i$, say $y_1$,  which
  implies that for $1 \leq i \leq 2r$, there exists a vertex $x_i\in
  X_1$ and a vertex $w_i\in Y_1$ so that $v_i$-$w_i$-$x_i$ form a path
  of length $2$.  Since $G$ is triangle-free, we have
  $w_i\not=w_{i+1}$ for all $i \in \{1,\ldots,2r-1\}$. Note that
  $x_i$-$w_i$-$v_i$-$v_{i+1}$-$w_{i+1}$-$x_{i+1}$ forms a path of
  length $5$ unless $x_i=x_{i+1}$.

Recall the admissible conditions: if  $\{u, v\} \subset X_i$ for
some $i \in \{1,2,3\}$, then $v\not\in N_G^1(u) \cup N_G^3(u) \cup
N_G^5(u)$. We must have $x_1=x_2=\cdots=x_{2r}$. Denote this common
vertex by $x$. Now have
$$\left|N_G^2(x)\cap B\right|\geq 2r.$$
Note $|N_G^2(x)|\leq 6$. We have $2r\leq 6$. The possible values for
$r$ are $2$ and $3$. If $r=2$, then $B$ is a $5$-cycle. Since $B$ is
in $T'$, we have $B\cap N^1_G(x)=\emptyset$; this is a contradiction
to Lemma \ref{c5}. If $r=3$, then $B$ is a $7$-cycle; this is a
contradiction to Lemma \ref{c7}.
\end{description}
The proof of  Lemma \ref{l:3colorable} is finished.
\hfill $\square$

\section{Partition into $42$ admissible sets}
\begin{theorem}\label{t5}
  Let $G$ be a triangle-free graph with maximum degree at most $3$.
  If $G$ is $2$-connected and $\textrm{girth}(G)\leq 6$,  then $G$ can be
  partitioned into at most $42$ admissible sets.
\end{theorem}
\noindent {\bf Proof of Theorem \ref{t5}:} 
We will define a proper coloring
$c\colon V(G) \to \{1,2,\ldots,126\}$ 
such that for $1 \leq i \leq 42$, the
$i$-th admissible set is $c^{-1}(\{3i-2,3i-1, 3i\})$. We refer to
$\{3i-2,3i-1, 3i\}$ as a color block for all $i \in
\{1,\ldots,42\}$.  Since $4\leq \ \textrm{girth(G)} \ \leq 6$, there
is a cycle $C$ of length $4$, $5$, or $6$. Let $v_{n-1}$ and $v_{n}$
be a pair of adjacent vertices of $C$. 
We assume that  $G\setminus\{v_{n-1},v_n\}$ is connected.
(If not, we start the greedy algorithm below from each of the components
of $G\setminus \{v_{n-1},v_n\}$ separately.)
 We can find a
vertex $v_1$ other than $v_{n-1}$ and $v_n$
 such that $G\setminus v_1$
is connected. Inductively, for each $i \in \{2,\ldots,n-2\}$, we can find a vertex $v_j$ other than $v_{n-1}$ and $v_n$ such that $G \setminus \{v_1,\ldots,v_{j-1}\}$ is connected. Therefore, we get an order of vertices
$v_1,v_2,\ldots, v_{n-1},v_n$ such that for $j=1,2,\ldots, n-2$,
the induced graph on $v_{j},\ldots, v_{n}$ is connected.

We color the vertices greedily. Assume we have colored $v_1,
v_2,\ldots, v_j$. For $v_{j+1}$, choose a color $h$ satisfying the
following:
\begin{enumerate}
\item For each $u\in N_G^1(v_{j+1}) \cap \{v_1, v_2,\ldots, v_j\}$,
we have  $h$ is not in the same block of $c(u)$.
\item For each $u\in \left(N_G^3(v_{j+1})\cup N_G^5(v_{j+1})\right)
\cap \{v_1, v_2,\ldots, v_j\}$, we have $h\not=c(u)$.
\item For each $u\in  \left(N_G^2(v_{j+1})\cup N_G^4(v_{j+1})\right)
\cap \{v_1, v_2,\ldots, v_j\}$, we have $h$ could equal  $c(u)$ but
not equal  the other two colors in the color block of $c(u)$.
\end{enumerate}
For $j\leq n-2$, there  is   at least one vertex in $N^1_G(v_{j+1})$
and one vertex in $N^2_G(v_{j+1})$ still uncolored. Thus
$|N_G^1(v_{j+1}) \cap \{v_1, v_2,\ldots,  v_j\}|\leq 2$, \hfill \\
$\left|N_G^2(v_{j+1}) \cap \{v_1, v_2,\ldots, v_j\}\right|\leq 5$,
$|N_G^3(v_{j+1})|\leq 12$, $|N_G^4(v_{j+1})|\leq 24$, and
$|N_G^5(v_{j+1})|\leq 48$. Since
$$3\times 2 + 2\times(5+24) + (12+48)=124<126,$$
it is always possible to color the vertex $v_{j+1}$ properly.

It remains to color $v_{n-1}$ and $v_n$ properly. Note both $v_{n-1}$
and $v_n$ are on the cycle $C$. Let us count color
redundancy according to the type of the cycle $C$.

\begin{description}
  \item[Case $C_4$:] For any vertex $v$ on $C_4$, there are two vertices
in $N_G^1(v)\cap N_G^3(v)$. We have $|N_G^2(v)|\leq 5$ and
$|N_G^4(v)|\leq 23$. Thus the number of colors  forbidden to be
assigned to  $v$ is at most
$$ 3\times 3 + 2\times(5+23) + (12+48) - 2=123<126.$$
 \item[Case $C_5$:] For any vertex $v$ on $C_5$, there are two vertices
in $N_G^1(v)\cap N_G^4(v)$. We also have $|N_G^5(v)|\leq 47$. The
number of  colors forbidden to be assigned to  $v$ is at most
$$ 3\times 3 + 2\times(6+24) + (12+47) - 2\times 2 =124<126.$$
\item[Case $C_6$:] For any vertex $v$ on $C_6$,
there are two vertices in $N_G^1(v)\cap N_G^5(v)$ and two vertices
in $N_G^2(v)\cap N_G^4(v)$. We  have $|N_G^3(v)|\leq 11$. Thus the
number of  colors  forbidden to be assigned to  $v$ is at most
$$ 3\times 3 + 2\times(6+24) + (11+48) - 2 - 2\times 2=122<126.$$
\end{description}
In each subcase, we can find a color for $v_{n-1}$ and $v_n$.

The 42 admissible sets can be obtained from the coloring $c$ as follows.
For $1\leq j\leq 42$, the $j$-th admissible set has
the following partition 
$$c^{-1}(3i-2) \cup c^{-1}(3i-1) \cup c^{-1}(3i).$$
Those are admissible sets by the construction of the coloring $c$.
\hfill $\square$

\noindent {\bf Proof of Theorem \ref{t1}:} Suppose that there exists
a graph $G$ which is triangle-free, $\Delta\leq3$, and
$\chi_f(G)>3-\frac{3}{43}$. Without loss of generality, we can
assume $G$ has the smallest number of edges
among all such graphs.
 Thus $G$ is $2$-connected and fractionally-critical.
If  $\textrm{girth}(G) \geq 7$,  Hatami and Zhu \cite{hz} showed
$\chi_f(G) \leq 2.78571 \leq 3 - \frac{3}{43}$. Contradiction!

If $\textrm{girth}(G) \leq 6$, Theorem \ref{t5} states that
$G$ can be partitioned into $42$ admissible sets. By Theorem \ref{t2},
we have $\chi_f(G)\leq 3-\frac{3}{k+1}=3-\frac{3}{43}$. Contradiction!
\hfill $\square$

{\bf Acknowledgment} We thank C.~C.~Heckman and anonymous referees
for very useful comments during the revision of this paper.

\end{document}